\documentclass[12pt]{amsart}

\input{xy}
\xyoption{all}
\usepackage{graphicx}
\usepackage{color}
\usepackage{verbatim}

\newtheorem{theorem}{Theorem} [section]
\newtheorem{prop}[theorem]{Proposition}
\newtheorem{lemma}[theorem]{Lemma}

\theoremstyle{definition}

\newtheorem{example}[theorem]{Example}

\theoremstyle{remark}
\newtheorem*{remark}{Remark}

\numberwithin{equation}{section}
\numberwithin{figure}{section}

\newcommand\C{{\mathbb C}}
\newcommand\CC{{\mathbb C}}

\renewcommand\P{{\mathbb P}}
\newcommand\PP{{\mathbb P}}
\newcommand\R{{\mathbb R}}
\newcommand\RR{{\mathbb R}}

\newcommand\Q{{\mathbb Q}}

\renewcommand\phi{\varphi}

\newcommand\cM{\mathcal{M}}

\newcommand\Gal{\operatorname{Gal}}

\renewcommand\Im {\operatorname{Im}}

\renewcommand\gcd {\operatorname{gcd}} 
\newcommand\Tor {\mathrm{Tor}} 
\newcommand\capacity {\operatorname{Cap}} 
\newcommand\Res {\operatorname{Res}} 

\newcommand\kbar {\overline{k}}
\newcommand\Qbar {\overline{\Q}}
\newcommand\kvbar {\overline{k}_v}
\newcommand\ksep {k^{\rm sep}}
\newcommand\Berk {\operatorname{Berk}}

\newcommand\hhat {\hat{h}}

\begin{document}

\title{Torsion points and the Latt\`es family}

\author{Laura De Marco}
\address{Department of Mathematics, Northwestern University, USA}
\email{demarco@math.northwestern.edu}
\author{Xiaoguang Wang}
\address{Department of Mathematics, Zhejiang University, P.R.China} 
\email{wxg688@163.com}
\author{Hexi Ye}
\address{Department of Mathematics, University of British Columbia, Canada}
\email{yehexi@math.ubc.ca}


\date{\today}

\begin{abstract}
We give a dynamical proof of a result of Masser and Zannier \cite{Masser:Zannier, Masser:Zannier:2}:  for any $a\not=b \in \Qbar\setminus\{0,1\}$, there are only finitely many parameters $t\in\C$ for which points $P_a = (a, \sqrt{a(a-1)(a-t)})$ and $P_b=(b,\sqrt{b(b-1)(b-t)})$ are both torsion on the Legendre elliptic curve $E_t = \{y^2 = x(x-1)(x-t)\}$.  Our method also gives the finiteness of parameters $t$ where both $P_a$ and $P_b$ have small N\'eron-Tate height.   A key ingredient in the proof is the equidistribution theorem of \cite{Baker:Rumely:equidistribution, FRL:equidistribution, ChambertLoir}.  For this, we prove two statements about the degree-4 Latt\`es family $f_t$ on $\P^1$:  (1) for each $c\in\C(t)$, the bifurcation measure $\mu_c$ for the pair $(f_t, c)$ has continuous potential across the singular parameters $t=0,1,\infty$; and (2) for distinct points $a,b\in\C\setminus\{0,1\}$, the bifurcation measures $\mu_a$ and $\mu_b$ cannot coincide.  Combining our methods with the result of \cite{Masser:Zannier:2}, we extend their conclusion to points $t$ of small height also for $a, b \in \C(t)$.  
\end{abstract}


\thanks{The research was supported by the National Science Foundation.}

\maketitle

\thispagestyle{empty}

\section{Introduction}

Masser and Zannier proved the following result about the Legendre family of elliptic curves $E_t$, defined by $\{y^2 = x(x-1)(x-t)\}$ for $t\in \C\setminus\{0,1\}$   \cite{Masser:Zannier:CR, Masser:Zannier, Masser:Zannier:2}.

\begin{theorem}  [Masser-Zannier]  \label{MZ}
Fix $a\not=b$ in $\Qbar\setminus\{0,1\}$.  There are only finitely many parameters $t\in\C$ for which the points $P_a(t) = (a,\sqrt{a(a-1)(a-t)})$ and $P_b(t) = (b,\sqrt{b(b-1)(b-t)})$ both have finite order on $E_t$.
\end{theorem}

\begin{remark}  Theorem \ref{MZ} is proved for $a=2$ and $b=3$ in \cite{Masser:Zannier:CR, Masser:Zannier}.  In \cite{Masser:Zannier:2}, the authors show that their result holds more generally.   Two given points $P_t$ and $Q_t$ (with $x$-coordinates in an algebraic closure of $\C(t)$) are simultaneously torsion for infinitely many parameters $t$ if and only if there are integers $m, n$ (not both zero) such that $[m]\cdot P_t = [n]\cdot Q_t$ for all $t$.  The motivation for this result is explained at length in these articles and in \cite{Zannier:book}.  Note that the point $P_a(t)$ with $a=0$ or $1$ will be torsion for all $t$, so these points are excluded in Theorem \ref{MZ}.
\end{remark}

Via the projection $E_t \to \P^1$ sending $(x,y)$ to $x$, Theorem \ref{MZ} can be reformulated in dynamical terms, about the family
\begin{equation} \label{Lattes family}
	f_t(z) = \frac{(z^2-t)^2}{4z(z-1)(z-t)}
\end{equation}
of complex rational functions on $\P^1$.  The map $f_t$ is induced by the multiplication-by-2 map on $E_t$; so a point $P$ is torsion on $E_t$ if and only if its projection to $\P^1$ has finite forward orbit for $f_t$.  As Zannier pointed out in \cite{Zannier:book}, Theorem \ref{MZ} is equivalent to the following statement:  {\em Fix $a\not=b$ in $\Qbar\setminus\{0,1\}$.  There are only finitely many parameters $t$ for which both $a$ and $b$ are preperiodic for the map $f_t$.}

In this article, we give an alternate proof of Theorem \ref{MZ} using dynamical and potential-theoretic methods.  In fact, we obtain a stronger result with these methods, a statement regarding points of small height in the spirit of the Bogomolov Conjecture (see \cite{Ullmo:Bogomolov, Zhang:Bogomolov}) and related to Zhang's Conjecture \cite[\S4]{Zhang:ICM}.  Fix $c\in \Qbar(t)$ with $c\not= 0,1,t$.  For each $t\in \Qbar\setminus\{0,1\}$, let $P_c(t)$ be a point on the elliptic curve $E_t$ with $x$-coordinate equal to $c(t)$.  Let $\hat{h}_c(t) := \hat{h}_{E_t}(P_c(t))$ be the N\'eron-Tate height of $P_c(t)$.  We define
	$$\Tor(c) = \{t\in \Qbar\setminus\{0,1\}:  c(t)  \mbox{ is a torsion point on } E_t \}$$
so that $t\in \Tor(c) \iff \hat{h}_c(t)=0$.  (See \S\ref{NT sum} and Lemma \ref{neron-tate} for details about $\hat{h}_c$.)  For constant points, we prove: 

\begin{theorem}  \label{synchrony}
Fix $a, b\in \Qbar$,  with $a, b \not= 0, 1$.  The following are equivalent:
\begin{enumerate}	
\item $|\Tor(a) \cap \Tor(b)| = \infty$;
\item  $\Tor(a) = \Tor(b)$; 
\item there is an infinite sequence $\{t_n\}\subset\Qbar$ so that $\hat{h}_a(t_n)\to 0$ and $\hat{h}_b(t_n) \to 0$; and 
\item $a=b$.
\end{enumerate}
\end{theorem}

If we combine our dynamical methods with the results of Masser and Zannier in \cite{Masser:Zannier:2}, we obtain the following result for $a$ and $b$ that are rational in $t$ with complex coefficients.  For $a,b\in\Qbar(t)$, we use the height functions $\hhat_a$ and $\hhat_b$ defined above.  If either $a$ or $b$ is not in $\Qbar(t)$, we let $\hat{h}_a$ and $\hat{h}_b$ be the canonical height functions defined over a finitely generated extension $k/\Qbar$ containing the coefficients of $a$ and $b$.  See \S\ref{function field height} for details.

\begin{theorem}  \label{synchrony over C}
Fix $a, b\in \C(t)$,  with $a, b \not= 0, 1, t$.  The following are equivalent:
\begin{enumerate}	
\item $|\Tor(a) \cap \Tor(b)| = \infty$;
\item  $\Tor(a) = \Tor(b)$; 
\item there is an infinite sequence $\{t_n\}\subset\kbar$ so that $\hat{h}_a(t_n)\to 0$ and $\hat{h}_b(t_n) \to 0$; and 
\item there exist nonzero integers $m$ and $n$ so that $[m]\cdot P + [n]\cdot Q \equiv 0$ on $E$.
\end{enumerate} 
\end{theorem}

The strategy of proof for Theorems \ref{synchrony} and \ref{synchrony over C} follows the ideas in the recent articles \cite{BD:preperiodic, BD:polyPCF, Ghioca:Hsia:Tucker, GHT:rational} which were themselves inspired by Theorem \ref{MZ} and related results; a key ingredient is the ``arithmetic equidistribution theorem" associated to height functions on $\P^1(\Qbar)$ of \cite{Baker:Rumely:equidistribution, FRL:equidistribution, ChambertLoir}.  However, the families of maps in \cite{BD:preperiodic, BD:polyPCF, Ghioca:Hsia:Tucker, GHT:rational}  all share certain technical features (e.g., the compactness of the bifurcation locus of a marked point), allowing the authors to apply the arithmetic equidistribution theorems more easily.  The Latt\`es family of (\ref{Lattes family}) requires a different approach to show that the adelic measures in question have continuous potentials on $\P^1$.  Most of this article is devoted to the arguments showing that the hypotheses of the equidistribution theorem are satisfied.

\subsection{Dynamical and potential-theoretic statements}  \label{results}
The Latt\`es map $f_t$ is postcritically finite for all $t$, meaning that each of its critical points has finite forward orbit, with postcritical set equal to $\{0,1,t,\infty\}$.  We consider {\em marked points} $c$ that are rational functions of $t$ (or constant $\infty$), viewed as holomorphic maps from the parameter space $\C\setminus\{0,1\}$ to the dynamical space $\P^1$.  We begin with an observation that the marked points in the postcritical set are the only rational points that are preperiodic for all $t$.

\begin{prop}  \label{persistent}
A marked point $c \in \C(t)\cup\{\infty\}$ is persistently preperiodic for $f_t$ if and only if $c(t)$ is equal to one of $0, 1, \infty, t$.
\end{prop}

\begin{remark}  
In terms of the Legendre family of elliptic curves $E_t$, Proposition \ref{persistent} states that $P_c = (c, \sqrt{c(c-1)(c-t)})$ is torsion for {\em all} $t\in \C\setminus\{0,1\}$ if and only if $c = 0, 1, t, \infty$.  This could also be deduced from more traditional methods in the study of elliptic curves.  Proposition \ref{first iterate} provides a more precise statement, computing the degrees of $t\mapsto f_t^n(c(t))$ for all $n\geq 1$ and all $c\in \C(t)$.  
\end{remark}

We will work in homogeneous coordinates for both the dynamical and parameter spaces.  We fix a homogenization of $f_t$,
	$$F_{t_1,t_2}(z,w) = ((t_1 w^2 - t_2 z^2)^2, 4t_2zw(w-z)(t_1w-t_2z)),$$		
for $(t_1, t_2)$ and $(z,w)$ in $\C^2$.

A {\em homogeneous lift} of a marked point $c$ means a pair of homogeneous polynomials $C = (c_1(t_1,t_2), c_2(t_1,t_2))$ (with no common factors) such that $c(t) = c_1(t,1)/c_2(t,1)$ for all $t$.  Given a pair of homogenous polynomials $F = (P,Q)$ of the same degree in $(t_1, t_2)$, we write $\deg F$ for their common degree, $\gcd(F)$ for the $\gcd$ of $P$ and $Q$ as polynomials in $(t_1, t_2)$, and $\Res(F)$ for the resultant of the pair $(P,Q)$.  We use the norm $\|(z,w)\| = \max\{|z|, |w|\}$ on $\C^2$.

\begin{theorem} \label{convergence}
Let $c\in \C(t)$ be a marked point $\not= 0,1,t$, and let $C$ be a homogeneous lift of $c$. We set $F_1 = F_{t_1,t_2}(C)/\gcd(F_{t_1,t_2}(C))$, $d = \deg F_1$, and $F_{n+1}=F_{t_1,t_2}(F_n)/t_2^2$ for all $n\geq 1$ so that ${\rm deg}(F_n)=4^{n-1} d$. Then the sequence
	$$\frac{1}{{\rm deg}(F_n)}\log\|F_n\|$$
converges locally uniformly for $(t_1,t_2) \in \mathbb{C}^2\setminus \{(0,0)\}$ to a continuous, plurisubharmonic function $G_C$, as $n\to \infty$.
\end{theorem}

The function $G_C$ of Theorem \ref{convergence} is a potential function for a probability measure on $\P^1$, the  {\em bifurcation measure} $\mu_c := \pi_* dd^c G_C$ associated to the marked point $c$, where $\pi(t_1,t_2) = t_1/t_2$ is the projection from $\C^2\setminus\{(0,0)\}$ to $\P^1$.   An important feature of $G_C$ is its continuity over parameters $t = 0,1,\infty$ where the family $f_t$ is not well defined.

\begin{theorem} \label{capacity}
Let $c$, $C$, and $F_1$ be as in Theorem \ref{convergence}.   The homogeneous capacity of the compact set
	$$K_C = \{(t_1,t_2)\in \mathbb{C}^2: G_C(t_1,t_2) \leq 0\}$$
is equal to
	$$\capacity(K_C)  = 4^{-\frac{1}{3 d}}\Big|\frac{(Q(1,1)-P(1,1))P(0,1)}{P(1,0)^2}\Big|^{-\frac{1}{3d^2}} |\Res(F_1)|^{-\frac{1}{d^2}}$$
where $F_1 = (P, Q)$.
\end{theorem}

\noindent
The homogeneous capacity in $\C^2$ was defined in \cite{D:lyap} and shown to coincide with the (square of the) transfinite diameter in \cite{Baker:Rumely:equidistribution}.  The proof of Theorem \ref{capacity} uses a relation between the capacity of $K_C$ and the resultants of $F_n$; compare \cite{DR:transfinite}.

We also prove non-archimedean analogs of Theorems \ref{convergence} and \ref{capacity}.  Specifically, when the marked point $c$ lies in $k(t)$ for a number field $k$,
we construct an {\em adelic} bifurcation measure $\mu = \{\mu_{c,v}\}$ where $v$ ranges over all places of $k$ and associated {\em canonical height function} $\hat{h}_\mu : \P^1(\bar{k}) \to \R_{\geq 0}$; we show that $\hat{h}_\mu$ agrees with the N\'eron-Tate height $\hat{h}_c$, up to a constant multiple (Proposition \ref{same heights}).  We are thus able to apply the arithmetic equidistribution theorems on $\P^1$  (as formulated in \cite{FRL:equidistribution, BRbook}) to study the set of torsion parameters $\Tor(c)$.

\begin{theorem}  \label{equidistribution}
Let $k$ be a number field.  For any $c \in k(t)$ ($\not= 0,1,t$), let $\{S_n\}$ be a non-repeating sequence of finite, $\Gal(\kbar/k)$-invariant sets in $\P^1(\kbar)$ for which
	$$\hat{h}_c(S_n) = \frac{1}{|S_n|} \sum_{t\, \in\, S_n} \hat{h}_c(t) \to 0.$$
Then the sets $S_n$ are equidistributed with respect to the bifurcation meausre $\mu_c$ on the parameter space $\C\setminus\{0,1\}$.  More precisely, the discrete probability measures
	$$\mu_n \; = \; \frac{1}{|S_n|} \; \sum_{t \,\in \,S_n} \; \delta_t$$
converge weakly to the measure $\mu_c$ as $n\to\infty$.  In particular, the torsion parameters $\Tor(c)$ are equidistributed with respect to $\mu_c$.  Moreover, for all places $v$ of $k$, the measures $\mu_n$ converge weakly to the bifurcation measure $\mu_{c,v}$ on the Berkovich space $\P^1_{\Berk,v}$.  
\end{theorem}

\begin{remark}
A version of Theorem \ref{equidistribution} holds for any $c\in\C(t)$, not only those points defined over a number field, and it is used to prove Theorem \ref{synchrony over C}.  However, we cannot directly equate the bifurcation measure $\mu_c$ on $\P^1(\C)$ with a component of the adelic measure $\{\mu_{c,v}\}$.  An additional argument would be needed to show that $\Tor(c)$ is equidistributed with respect to $\mu_c$; compare \cite[Theorem 1.5]{Yuan:Zhang:II}.  
\end{remark}

The measure of maximal entropy for the Latt\`es map $f_t$ is equal to the Haar measure on the torus, projected to $\P^1$; consequently, we are able to provide an explicit integral expression for the potential function $G_C$ of Theorem \ref{convergence}.  We use the expression for the density function $\rho_\Lambda(t)$ for the hyperbolic metric in $\Lambda=\C\setminus\{0,1\}$ provided in \cite[Theorem 4.13]{McMullen:course}.

\begin{theorem} \label{integral formula}
The function $G_C$ of Theorem \ref{convergence} has the following explicit expression
\begin{eqnarray*}
 G_C(t_1,t_2)&=&\frac{2}{d}\left(|t_2|D(t_1/t_2) \cdot \int_{\mathbb{P}^1}\frac{\log|c_1-c_2\zeta|}{|\zeta(\zeta-1)(t_2\zeta-t_1)|}|d\zeta|^2+\log|t_2|\right)\\
&& \qquad -\frac{1}{d}\log|{\rm gcd}(F_{t_1,t_2}(C))|
 \end{eqnarray*}
where $C=(c_1,c_2)$ and $D(t)=4|t(t-1)|\rho_\Lambda(t)$.
\end{theorem}

As a consequence of Theorem \ref{integral formula}, and towards the proof of Theorem \ref{MZ}, we show:

\begin{prop}  \label{distinct measures}
For constant marked points $a, b \in \C\setminus\{0,1\}$, the bifurcation measures satisfy $\mu_a = \mu_b$ if and only if $a=b$.
\end{prop}

\subsection{Proof of Theorem \ref{MZ}}
Fix $a, b \in \Qbar\setminus\{0,1\}$.  Let $k$ be a number field containing $a$ and $b$. Suppose $\{t_n\} \subset \C$ is an infinite sequence of parameters at which both $a$ and $b$ are preperiodic for $f_t$.  Then $t_n\in \kbar$ and $\hat{h}_a(t_n) = \hat{h}_b(t_n) = 0$ for all $n$.  By Theorem \ref{equidistribution}, the $\Gal(\kbar/k)$-orbits of these parameters are equidistributed with respect to both $\mu_a$ and $\mu_b$.  Consequently $\mu_a = \mu_b$, so by Proposition \ref{distinct measures} we must have $a=b$.  This completes the proof of Theorem \ref{MZ}.

The proof of Theorem \ref{synchrony} is given in \S\ref{proof of synchrony}, and the proof of Theorem \ref{synchrony over C} is completed in \S\ref{complex case}.

\subsection{Acknowledgements}
We are greatly indebted to Matthew Baker, who introduced us to this problem; he has been generous with his mathematical ideas.  We would also like to thank Dragos Ghioca and Tom Tucker and the anonymous referees for many helpful comments and suggestions.  Finally, we thank Umberto Zannier for the encouragement to develop this alternate proof of his theorem with David Masser.

\bigskip
\section{Resultants and Capacity}

For the analysis of the bifurcation measures on $\P^1$, we will typically work in homogeneous coordinates on $\C^2$.   In this section, we provide background on the resultant of two homogeneous polynomials (in two variables) and the homogeneous capacity for compact sets in $\C^2$.

\subsection{The resultant}
Let $F:\mathbb{C}^2\rightarrow \mathbb{C}^2$ be a homogeneous polynomial map of degree $d\geq 1$. It can be written as $F(z_1,z_2)=(a_0z^d_1+\cdots+ a_dz^d_2, b_0z^d_1+\cdots+ b_dz^d_2)$.  The resultant of $F$ is the resultant of the pair of polynomials defining $F$, computed explicitly as the determinant of the $(2d)\times (2d)$ matrix,
  $${\rm Res}(F)=\left|
\begin{array} {cccccccc}
a_0 & a_1& \cdots & a_{d-1} & a_d &0 &\cdots & 0\\
 0 & a_0 & a_1& \cdots & a_{d-1} & a_d & \cdots &0\\
  & & & \vdots & & & & \\
   0 & 0& \cdots & a_0 & a_1& \cdots&  a_{d-1} & a_d\\
  b_0 & b_1&\cdots & b_{d-1} & b_d &0 &\cdots & 0\\
 0 & b_0 & b_1& \cdots & b_{d-1} & b_d & \cdots &0\\
  & & & \vdots & & & & \\
   0 & 0& \cdots & b_0 & b_1&\cdots&  b_{d-1} & b_d
\end{array}
\right|.$$

Given a pair of points $z=(z_1,z_2)$ and $w=(w_1,w_2)$
in $\mathbb{C}^2$, we may defined the wedge product by $z\wedge w=z_1w_2-z_2w_1$. Then  $F$ can be
factored as
$$F(z)=\Big(\prod z\wedge \alpha_i, \prod z\wedge \beta_j\Big)$$
 for some points $\alpha_i,\beta_j\in\mathbb{C}^2$, and the resultant satisfies
 $${\rm Res}(F)=\prod_{i,j}\alpha_i\wedge\beta_j.$$

Note that ${\rm Res}(F)\neq0$ iff $F^{-1}(0,0)=\{(0,0)\}$. In this case, we say $F$ is non-degenerate.

\subsection{The homogeneous capacity}
We will consider compact sets $K\subset \C^2$ that are circled and pseudoconvex:  these are sets of the form
	$$K = \{(z,w)\in \mathbb{C}^2: G_K(z,w) \leq 0\}$$
for continuous, plurisubharmonic functions $G_K :  \C^2\setminus\{(0,0)\}  \to \R$ such that
	$$G_K(\alpha z, \alpha w) = G_K(z, w) + \log|\alpha|$$
for all $\alpha\in\C^*$; see \cite[\S3]{D:lyap}.  Such functions are (homogeneous) potential functions for probability measures on $\P^1$ \cite[Theorem 5.9]{Fornaess:Sibony}.

Set $G_K^+=\max\{G_K,0\}$.  The Levi measure of $K$ is defined by
	$$\mu_K=dd^cG_K^+\wedge dd^cG_K^+.$$
It is known that $\mu_K$ is a probability measure supported on $\partial K=\{G=0\}$. The homogeneous capacity of $K$ is defined by
$${\rm Cap}(K)=\exp\Big(\iint\log|\zeta\wedge \xi| d\mu_K(\zeta)d\mu_K(\xi)\Big).$$
This capacity was introduced in \cite{D:lyap} and shown in \cite{Baker:Rumely:equidistribution} to satisfy $\capacity(K) = (d_\infty(K))^2$, where $d_\infty$ is the transfinite diameter in $\C^2$.   We will make use of the following pullback formula for the homogeneous capacity.

\begin{prop} \label{pull-back}
Let $K\subset\mathbb{C}^2$ be a  compact circled and pseudoconvex set, and $F$ be a non-degenerate homogenous polynomial map of
degree $d\geq 1$,  we have
$$\capacity(F^{-1}K)=|{\rm Res}(F)|^{-1/d^2}{\rm Cap}(K)^{1/d}$$
\end{prop}

\noindent
A more general version of Proposition \ref{pull-back} was proved in \cite{DR:transfinite}; see also \cite[Corollary 6.3]{Berman:Boucksom}.  In dimension 2 with homogeneous maps and circled sets, this pullback formula has a more elementary proof that we include here.

\begin{proof} Set $e_1=(0,1)$ and $e_2=(1,0)$. By \cite[Lemma 4.3]{D:lyap},
$$|{\rm Res}(F)|^{d^2}\prod_{F(z)=e_1}\prod_{F(w)=e_2}|z\wedge w|=1.$$
For any $\zeta,\xi\in\mathbb{C}^2$ with $\zeta\wedge\xi\neq0$, take a linear map $\ell:\mathbb{C}^2\rightarrow\mathbb{C}^2$ such that
$\ell(e_1)=\zeta$ and $\ell(e_2)=\xi$. In fact, $\ell(z_1,z_2)=(\xi_1z_1+\zeta_1z_2, \xi_2z_1+\zeta_2z_2)$. Then,
\begin{eqnarray*}
\prod_{F(z)=\zeta}\prod_{F(w)=\xi}|z\wedge w|&=&\prod_{\ell^{-1}F(z)=e_1}\prod_{\ell^{-1}F(w)=e_2}|z\wedge w|\\
&=&|{\rm Res}(\ell^{-1}F)|^{-d^2}\\
&=&|\det(\ell^{-1})^{d} {\rm Res}(F)|^{-d^2}\\
&=&|{\rm Res}(F)|^{-d^2}|\zeta\wedge\xi|^{d^3}.
\end{eqnarray*}

Let $G_K$ be the defining function for $K$, and set $G^+_{F^{-1}K}=\frac{1}{d}G_K^+\circ F$. The Levi measure of $F^{-1}K$ is
$\mu_{F^{-1}K}=dd^cG^+_{F^{-1}K}\wedge dd^cG^+_{F^{-1}K}$.  Let $\Delta$  be the diagonal of $\mathbb{C}^2$. We have
\begin{eqnarray*}
\log{\rm Cap}(K)&=&\iint\log|\zeta \wedge \xi| d\mu_K(\zeta)d\mu_K(\xi)\\
&=&\iint_{\mathbb{C}^2\setminus\Delta}\log|\zeta \wedge \xi| d\mu_K(\zeta)d\mu_K(\xi)\\
&=&\iint_{\mathbb{C}^2\setminus\Delta}\left(\frac{1}{d^3}\sum_{F(z)=\zeta}\sum_{F(w)=\xi}\log|z\wedge w|+\frac{1}{d}\log|{\rm Res}(F)|\right)d\mu_K(\zeta)d\mu_K(\xi)\\
&=&\frac{1}{d}\log|{\rm Res}(F)|+d\iint_{\mathbb{C}^2\setminus F^{-1}\Delta} \log|z\wedge w|d\mu_{F^{-1}K}(z)d\mu_{F^{-1}K}(w)\\
&=&\frac{1}{d}\log|{\rm Res}(F)|+d\log{\rm Cap}(F^{-1}K).
\end{eqnarray*}
Equivalently,
$${\rm Cap}(F^{-1}K)=|{\rm Res}(F)|^{-1/d^2}{\rm Cap}(K)^{1/d}.$$
\end{proof}

\subsection{Evaluating capacity via resultant}

Let $G_K, K$ be defined as in the previous section.

\begin{theorem}\label{cap=res} Suppose that there is a sequence of non-degenerate  homogenous polynomials $\{F_n:\mathbb{C}^2\rightarrow\mathbb{C}^2\}_{n\geq1}$ such that $$\frac{1}{d(F_n)}\log\|F_n\|$$
converges locally uniformly to $G_K$ in $\mathbb{C}^2-\{(0,0)\}$ as $n\rightarrow+\infty$, with $d(F_n)$ the degree of $F_n$. Then we have
$$ {\rm Cap}(K)=\lim_{n\rightarrow\infty} |{\rm Res}(F_n)|^{-1/d(F_n)^2}.$$
\end{theorem}

\noindent
Theorem \ref{cap=res} implies the existence of the limit $\lim_{n\rightarrow\infty} |{\rm Res}(F_n)|^{-1/d(F_n)^2}$.  Note that we do not require growth of the degree $d(F_n)$.

  \begin{proof} By assumption, the sequence $\frac{1}{d(F_n)}\log\|F_n\|$ converges uniformly to $G_K$ on $\partial K$, where $G_K$ vanishes. So for any $\delta>0$, there is an integer $n(\delta)$ such that for all $n\geq n(\delta)$ and  all $(z_1,z_2)\in \partial K$,
$$-\delta\leq \frac{1}{d(F_n)}\log\|F_{n}(z_1,z_2)\|\leq\delta.$$
The homogeneity of $F_n$ implies that for any $\alpha\neq0$, one has $F_n(\alpha z_1,\alpha z_2)=\alpha^{d(F_n)}F_{n}(z_1,z_2)$. It follows that
$$\frac{1}{d(F_n)}\log\|F_{n}(\alpha z_1,\alpha z_2)\|=\frac{1}{d(F_n)}\log\|F_{n}(z_1,z_2)\|+\log|\alpha|$$
 As a  consequence,  for any $n\geq n(\delta)$, we have
 \begin{equation*}
\frac{1}{d(F_n)}\log\|F_{n}(z_1,z_2)\| \; \begin{cases}
 \leq\delta,\ &\forall (z_1,z_2)\in K,\\
>-\delta,\ &\forall (z_1,z_2)\in \mathbb{C}^2 \setminus K.
\end{cases}
\end{equation*}
It follows that
$$F_{n}^{-1}({\mathbb{D}^2(e^{-d(F_n)\delta})})\subset K \subset F_n^{-1}({\mathbb{D}^2(e^{d(F_n)\delta})})$$
where $\mathbb{D}^2(r)$ is the  polydisk $\{(z_1,z_2); |z_1|\leq r, |z_2|\leq r\}$, whose capacity is  $r^2$ \cite[\S4]{D:lyap}.
 Thus by Proposition \ref{pull-back} and monotonicity of capacity,
 $$ e^{-2\delta} |{\rm Res}(F_n)|^{-1/d(F_n)^2}\leq {\rm Cap}(K)\leq e^{2\delta} |{\rm Res}(F_n)|^{-1/d(F_n)^2}, \ \ \forall n\geq n(\delta).$$
Since $\delta>0$ can be arbitrarily small, we have
$$\limsup_{n\rightarrow\infty} |{\rm Res}(F_n)|^{-1/d(F_n)^2}\leq {\rm Cap}(K)\leq \liminf_{n\rightarrow\infty} |{\rm Res}(F_n)|^{-1/d(F_n)^2}.$$
The conclusion then follows.
\end{proof}

\bigskip
\section{Marked points and capacity}
\label{growth}

In this section, we provide a proof of Proposition \ref{persistent}; a more precise version is stated as Proposition \ref{first iterate} below.  Also, assuming Theorem \ref{convergence}, we provide the proof of Theorem \ref{capacity}.

\subsection{Degree growth}
Recall that we have
	$$f_t(z) = \frac{(z^2-t)^2}{4z(z-1)(z-t)},$$
for all $t\in\C\setminus\{0,1\}$, and a homogenization
	$$F_{t_1,t_2}(z,w) =  ((t_1 w^2 - t_2 z^2)^2, 4t_2zw(w-z)(t_1w-t_2z))$$
for all $(t_1, t_2) \in \C^2$ such that $t = t_1/t_2 \not= 0, 1, \infty$.   A marked point is an element $c$ of $\C(t)\cup \{\infty\}$, viewed as a holomorphic map from the parameter space $t\in\C\setminus\{0,1\}$ to $\P^1$.

\begin{prop} \label{first iterate}
Let $C$ be a homogenous lift of any marked point not equal to $0, 1, t, \infty$.  We set $F_1 = F(C)/\gcd(F(C)) =: (P_1, Q_1)$.  Then $d := \deg F_1 \geq 2$, $P_1(0,1) \not=0$, $P_1(1,0) \not=0$, $P_1(1,1) \not= Q_1(1,1)$, and $Q_1(1,0) = 0$.   Furthermore, for all $n\geq 1$, the map
	$$F_{n+1}=F_{t_1,t_2}(F_n)/t_2^2$$
has nonzero resultant and degree ${\rm deg}(F_{n+1})=4^n d$.
\end{prop}

\proof
First, if $C = (z,w)$ is constant (and $\not= (0,\alpha),(\alpha,0), (\alpha,\alpha)$ for any $\alpha\in\C$), then
	$$F_1(t_1, t_2) =   ((t_1 w^2 - t_2 z^2)^2, 4t_2zw(w-z)(t_1w-t_2z)),$$
which has degree 2 and nonzero resultant.  The conditions on $P_1$ and $Q_1$ are immediate.

If $C = (A,B)$ has degree $m \geq 1$ with no common factors, we consider multiple cases.  Write
	$$F(C) = ((t_1 B^2 - t_2 A^2)^2, 4t_2AB(B-A)(t_1B-t_2A)).$$
If $(A,B) = (w t_1, z t_2)$, with $(z,w) \not= (0,\alpha),(\alpha,0), (\alpha,\alpha)$, then a computation shows that $\gcd(F(C)) = t_1^2t_2^2$ and
	$$F_1 = F_{t_1, t_2}(z,w).$$
We have already seen that this $F_1$ has degree 2 with the desired conditions on $P_1$ and $Q_1$.

Still assuming $m\geq 1$, write $A = t_1^l A'$ and $B = t_2^k B'$ for $k,l\geq 0$ maximal.  Then
	$$t_1^{\min\{2, 2l\}} t_2^{\min\{2, 2k\}} | (t_1 B^2 - t_2 A^2)^2$$
and these are the maximal such powers of $t_1$ and $t_2$; furthermore, $(t_1 B^2 - t_2 A^2)^2$ can share no other factors with $A$ or $B$.  Note also that
	$$t_1^{\min\{1, l\}} t_2^{\min\{1, k\}} | (t_1B - t_2A) \quad\mbox{ and } \quad t_1^lt_2^{k+1}|(t_2AB).$$
For $k=0$, it follows immediately that $\deg(F_1) \geq \deg(t_2B) =m+1 \geq 2$.  For $l=0$ and $k\geq 1$, we see that $\deg(F_1) \geq \deg(AB) = 2m \geq 2$.  For $k,l\geq 1$, we conclude that $\deg(F_1) \geq \deg(AB/t_1) \geq 2m-1$, which is $\geq 2$ unless $m=1$; but the case $l = k = m =1$ was treated above.

Note, in each case, that $t_1^{\min\{2, 2l\}} t_2^{\min\{2, 2k\}}$ divides the $\gcd$ of $F(C)$, and these are the maximal such powers of $t_1$ and $t_2$ to do so.  After cancellation, $t_2$ will still divide the second term in $F_1$, and we see that $Q_1(1,0) = 0$, that $P_1(1,0)\not=0$, and that $P_1(0,1) \not=0$.

It remains to show that $P_1(1,1) \not= Q_1(1,1)$.  Set $a = A(1,1)$ and $b = B(1,1)$.  Evaluating at $(t_1, t_2) = (1,1)$, we have
	$$F(C)(1,1) = ((b-a)^2(b+a)^2, 4ab(b-a)^2).$$
If $a\not=b$, then the two sides are clearly distinct, so also $P_1(1,1) \not= Q_1(1,1)$ in this case.  Finally, suppose $a=b$.  Note that $a \not=0$ because $A$ and $B$ have no common factors.  Then $(t_1 - t_2)^2$ must divide both $(B-A)(t_1B-t_2A)$, and $(t_1 B^2 - t_2 A^2)^2$.  We can compute $P_1(1,1)$ and $Q_1(1,1)$ by evaluating the limit of $F(C(t,1)) /(t-1)^2$ as $t\to 1$.  Expanding in a series around $t=1$, write $A(t,1) = a + a' (t-1) + \cdots$ and $B(t,1) = a + b'(t-1) + \cdots$.  Then $\lim_{t\to 1} (t B^2 - A^2)^2/(t-1)^2 = a^2(a - 2a' + 2b')^2$, and $\lim_{t\to 1} 4AB(B-A)(tB-A)/(t-1)^2 = -4a^2 (a'-b')(a - a'+b')$.  Since $a\not=0$, we see that $P_1(1,1) \not= Q_1(1,1)$.

Now consider $F_{n+1} = F_{t_1, t_2}(F_n)/t_2^2$ for all $n\geq 1$.  Applying the arguments above, $P_1(0,1) \not=0$ and $P_1(1,1) \not= Q_1(1,1)$ imply that neither $t_1$ nor $(t_1-t_2)$ may be factors of $\gcd(F(F_1))$, and $Q_1(1,0) = 0$ implies that $t_2^2$ is the highest power of $t_2$ that divides the gcd.  The fact that $F_{t_1, t_2}(z,w)$ is a non-degenerate homogeneous polynomial in $(z,w)$ for all $t_1/t_2 \not= 0,1, \infty$ implies that there can be no other factors of $\gcd(F(F_1))$.  Furthermore, the formula for $F_{t_1, t_2}$ shows immediately that $F_2$ must have degree $4\cdot \deg F_1$.  Therefore, $F_2 = (P_2, Q_2)$ satisfies the conditions of the proposition, and again we may conclude that $P_2(0,1), P_2(1,0) \not=0$, $Q_2(1,0)=0$, and $P_2(1,1) \not= Q_2(1,1)$.  Continuing inductively, the proposition is proved.
\qed

\subsection{Proof of Proposition \ref{persistent}}
Let $c$ be any marked point, not equal to $0, 1, t, \infty$.  Let $C$ be any homogeneous lift of $c$.  From Proposition \ref{first iterate}, the degree of $f^n_t(c(t))$ in $t$ is growing exponentially with $n$, with no cancellation in numerator and denominator after the first iterate.  Consequently, $c(t)$ cannot be preperiodic.

\subsection{Computing resultants}
The next proposition computes the resultants of the polynomial maps $F_n(t_1, t_2)$ appearing in Propostion \ref{first iterate}.

\begin{prop} \label{res-lat}
Let $G=(P,t_2Q)$ be a homogenous polynomial in $t=(t_1,t_2)$ of degree $d$.
Then the resultant of $F=F_{t_1,t_2}(G)/t_2^2$ is given by
$${\rm Res}(F)=4^{4d}[Q(1,1)-P(1,1)]^4 P(0,1)^4P(1,0)^{-8}{\rm Res}(G)^{16}.$$
In particular, the $F_n$ of Proposition \ref{first iterate} has resultant
	$$\Res(F_n) = 4^{A_nd}\Bigg(\frac{(P_{1}(1,1)-Q_{1}(1,1))P_1(0,1)}{P_1(1,0)^2}\Bigg)^{A_n}{\rm Res}(F_{1})^{4^{2(n-1)}},$$
for all $n\geq 1$, with $A_n=4^{n-1}+4^n+\cdots+4^{2n-3}=\frac{4^{2(n-1)}-4^{n-1}}{3}$.
\end{prop}

\proof
Let $R,S$ be two homogeneous polynomials in $t=(t_1,t_2)$, not necessarily having the same degree. Then each one can be factored as
$\prod t\wedge \zeta_i$. The resultant of the pair $(R,S)$ is equal to
$${\rm Res}(R,S)=\prod_{R(\xi)=0}\prod_{S(\zeta)=0}\xi\wedge\zeta=\prod_{R(\xi)=0}S(\xi)=(-1)^{{\rm deg}(R){\rm deg}(S)}\prod_{S(\zeta)=0}R(\zeta).$$
A basic fact is ${\rm Res}(R_1R_2,S)={\rm Res}(R_1,S){\rm Res}(R_2,S)$. By this fact, the resultant of $F=( (P^2 - t_1t_2Q^2)^2, 4t_2 PQ(P-t_2Q)(P-t_1Q))$
satisfies ${\rm Res}(F)=4^{4d}(I_1I_2I_3I_4)^2$, where $I_j$  are defined and evaluated as follows:
\begin{eqnarray*}
I_1&:=&{\rm Res}(t_2Q,  P^2-t_1t_2Q^2)=\prod_{t_2Q=0} (P^2-t_1t_2Q^2)\\
&=&\prod_{t_2Q=0}P^2={\rm Res}(G)^2,
\end{eqnarray*}
\begin{eqnarray*}
 I_2&:=&{\rm Res}(P,  P^2-t_1t_2Q^2)=\prod_{P=0} (P^2-t_1t_2Q^2)\\
&=&\prod_{P=0}(-t_1t_2Q^2)={\rm Res}(G)^2\prod_{P=0}t_1/\prod_{P=0}(-t_2)  \\
&=&{\rm Res}(G)^2 P(0,1)P(1,0)^{-1},
\end{eqnarray*}
\begin{eqnarray*}
 I_3&:=&{\rm Res}(P-t_2Q,  P^2-t_1t_2Q^2)=\prod_{P-t_2Q=0} (P^2-t_1t_2Q^2)\\
&=&\prod_{P-t_2Q=0}(t_2-t_1)t_2Q^2=\prod_{(t_2-t_1)t_2Q^2=0}(P-t_2Q) \\
&=&[P(1,1)-Q(1,1)]P(1,0)^{-1}\prod_{t_2^2Q^2=0}P\\
&=&[P(1,1)-Q(1,1)]P(1,0)^{-1}{\rm Res}(G)^2,
\end{eqnarray*}
\begin{eqnarray*}
 I_4&:=&{\rm Res}(P-t_1Q,  P^2-t_1t_2Q^2)=\prod_{P-t_1Q=0} (P^2-t_1t_2Q^2)\\
&=&\prod_{P-t_1Q=0}(t_1-t_2)t_1Q^2=\prod_{(t_1-t_2)t_1Q^2=0}(P-t_1Q)\\
&=&[P(1,1)-Q(1,1)]P(0,1)P(1,0)^{-2}{\rm Res}(G)^2.
\end{eqnarray*}
So we have
$${\rm Res}(F)=4^{4d}[P(1,1)-Q(1,1)]^4 P(0,1)^4P(1,0)^{-8}{\rm Res}(G)^{16}.$$

For $n\geq1$, let $F_n = (P_n, t_2Q_n)$ be the pair of homogeneous polynomials from Proposition \ref{first iterate}. By Proposition \ref{res-lat},
$${\rm Res}(F_{n+1})=\frac{4^{4{\rm deg}(F_{n})}[P_{n}(1,1)-Q_{n}(1,1)]^4P_{n}(0,1)^4{\rm Res}(F_{n})^{16}}{P_{n}(1,0)^8}.$$
Note that ${\rm deg}(F_{n})=4^{n-1}d$. By the relation
 $$F_{n+1}=(P_{n+1}, t_2Q_{n+1})=((P_n^2-t_1t_2Q_n^2)^2,4t_2 P_nQ_n(P_n-t_2Q_n)(P_n-t_1Q_n))$$
we get
$$P_{n}(1,0)=P_{n-1}(1,0)^4=P_{1}(1,0)^{4^{n-1}}, $$
$$P_{n}(0,1)=P_{n-1}(0,1)^4=P_{1}(0,1)^{4^{n-1}}, $$
and
$$P_{n}(1,1)-Q_{n}(1,1)=(P_{n-1}(1,1)-Q_{n-1}(1,1))^4=(P_{1}(1,1)-Q_{1}(1,1))^{4^{n-1}}.$$
So we have
\begin{eqnarray*}
{\rm Res}(F_{n}) &=& 4^{4^{n-1}d}\Bigg(\frac{(P_{1}(1,1)-Q_{1}(1,1))P_1(0,1)}{P_1(1,0)^2}\Bigg)^{4^{n-1}}{\rm Res}(F_{n-1})^{4^2} \\
 & =& 4^{A_nd}\Bigg(\frac{(P_{1}(1,1)-Q_{1}(1,1))P_1(0,1)}{P_1(1,0)^2}\Bigg)^{A_n}{\rm Res}(F_{1})^{4^{2(n-1)}},
\end{eqnarray*}
where
$$A_n=4^{n-1}+4^n+\cdots+4^{2n-3}=\frac{4^{2(n-1)}-4^{n-1}}{3}.$$
\qed

\begin{example}
Suppose the starting point $C = (z,w)$ is constant (with $z/w \not= 0, 1, \infty$).
 Then	
 $$\Res(F_n) = [2w(w-z)]^{a_n} z^{4^{2n}-2a_n}$$
where $a_n = 4^n + 4^{n+1} + \cdots + 4^{2n-1}$.
\end{example}

\subsection{Proof of Theorem \ref{capacity}}
Assuming Theorem \ref{convergence}, the proof is immediate from Proposition \ref{res-lat} and Theorem \ref{cap=res}.

\bigskip
\section{Convergence proof}
\label{estimates}

This section is devoted to the proof of Theorem \ref{convergence}.  The proof is lengthy, though each step is elementary.  We include details, especially so that we may compare this proof to the non-archimedean version needed for our proof of Theorem \ref{equidistribution}.

Recall that
	$$F_{t_1,t_2}(z,w) =  ((t_1 w^2 - t_2 z^2)^2, 4t_2zw(w-z)(t_1w-t_2z))$$
for all $(t_1, t_2), (z,w) \in \C^2$.  For starting point $C(t_1, t_2)$, a homogeneous lift of a marked point $c(t) \not= 0, 1, t, \infty$, we defined
	$$F_1 = F_{t_1,t_2}(C)/\gcd(F_{t_1,t_2}(C)) =: (P_1(t_1, t_2), t_2Q_1(t_1, t_2))$$
and
	$$F_{n+1}=F_{t_1,t_2}(F_n)/t_2^2$$
for all $n\geq 1$.  We set $d = \deg F_1$ so that ${\rm deg}(F_n)=4^{n-1} d$.  Our goal is to prove that the sequence
	$$\frac{1}{d\cdot 4^{n-1}}\log\|F_n\|$$
converges locally uniformly for $(t_1,t_2) \in \mathbb{C}^2\setminus \{(0,0)\}$ to a continuous, plurisubharmonic function $G_C$, as $n\to \infty$.  It will be convenient to work with the norm
	$$\|(z,w)\| := \max\{|z|, |w|\}$$
for $(z,w) \in \C^2$.

The proof of Theorem \ref{convergence} also shows:

\begin{prop} \label{values of G_C}
Under the hypotheses of Theorem \ref{convergence}, the function $G_C$ over the degenerate parameters $t=t_1/t_2=0,1,\infty$ satisfies
$$G_C(1,0) = \frac{1}{d} \log |P_1(1,0)|, \qquad G_C(0,1) = \frac{1}{d} \log|P_1(0,1)|, $$
and
$$G_C(1,1) = \frac{1}{d} \log |P_1(1,1)-Q_1(1,1)|.$$
\end{prop}

\subsection{Convergence for non-degenerate parameters}
As a map in $(z,w)$, $F_{t_1, t_2}$ has a nonzero resultant if and only if $t_1/t_2 \not= 0, 1, \infty$.  Therefore, applying the (now-standard) methods from complex dynamics, we know that the sequence
	 $$\frac{1}{d \cdot 4^{n-1}}\log\|F_n\|$$
converges locally uniformly in $\C^2\setminus \{(t_1,t_2): t_1t_2(t_1-t_2)=0\}$.  See \cite{Hubbard:Papadopol} or \cite{Fornaess:Sibony} for details.

As the maps $F_n$ are homogeneous, the proof of Theorem \ref{convergence} reduces to showing that the sequences
$$\frac{1}{d\cdot 4^{n-1}}\log\|F_n(1,t)\|, \quad  \frac{1}{d \cdot 4^{n-1}}\log\|F_n(t,1)\|, \quad \mbox{ and } \quad \frac{1}{d\cdot 4^{n-1}} \log\|F_n(1, 1-t)\|$$
converge locally uniformly in a small neighborhood of $t=0$.

\subsection{Convergence near $t_2=0$}
\label{near infinity}
We will prove that for any sufficiently small $\epsilon>0$, there are small $\delta_{\epsilon}>0$ and  integer $N_{\epsilon}$ such that for any $n\geq N_{\epsilon}$ and $|t|<\delta_{\epsilon}$ we have
\begin{equation}\label{lattes Fn 1,0}-\epsilon<\frac{1}{{\rm deg}(F_n)}\log\|F_n(1, t)\|-\frac{\log |P_1(1,0)|}{d}<\epsilon.
\end{equation}
These inequalities will be guaranteed by the following lemmas.

We first study the coefficients of $F_n(1,t)$. Write
   $$F_n(1,t)=(P_n(1,t), tQ_n(1,t))=( A_n+O(t), t(B_n+O(t))).$$

\begin{lemma}\label{t=0 coefficients growth}
For all $n\geq 1$, we have $A_n=P_1(1,0)^{4^{n-1}}$.  For the sequence $B_n$,
   $$\lim_{n\to \infty} \sup |B_n|^{1/4^{n-1}}\leq |P_1(1,0)|.$$
\end{lemma}

\proof From the inductive formula $F_{n+1}(1,t)=F_{1,t}(F_n(1,t))/t^2$,
we see that $A_{n+1}=A_n^4$ and $A_1=P_1(1,0)$. So for any $n\geq 1$, $A_n=P_1(1,0)^{4^{n-1}}$.

Let $b_1=0$ and $b_n=2b_{n-1}+1$.  Then inductively, we can write
       $$B_n=P_1(1,0)^{4^{n-1}-b_n-1}B_n^*,$$
with $B_1^*=Q_1(1,0)$ and $B_n^*=4B_{n-1}^*(P_1(1,0)^{b_{n-1}+1}-B_{n-1}^*)$. As $B_n^*$ grows quadratically when $n$ increases by one, we obtain $\limsup_{n\to \infty}  |B_n^*|^{1/4^{n-1}}\leq 1$. Consequently, we have $\limsup_{n\to \infty} |B_n|^{1/4^{n-1}}\leq |P_1(1,0)|.$\qed

\medskip
The growth of the coefficients of $A_n$ and $B_n$ in Lemma \ref{t=0 coefficients growth} provides a uniform upper bound on the size of $\frac{1}{{\rm deg}(F_n)}\log\|F_n(1,t)\|$ for small $t$:

\begin{lemma}\label{Lattes t=0 upper bound}
For any given $\epsilon>0$, there exists a $\delta>0$ and an integer $N>0$ so that
   $$\frac{1}{{\rm deg}(F_n)}\log\|F_n(1,t)\|-\frac{\log |P_1(1,0)|}{d}<\epsilon$$
for all $|t|<\delta$ and all $n\geq N$.
\end{lemma}

\proof We define the polynomials $p_n(t)$ and $q_n(t)$ by
     $$F_n(1,t)=:(p_n(t), tq_n(t))$$
so that $p_n(0)=A_n$ and $q_n(0)=B_n$. By Lemma \ref{t=0 coefficients growth}, there is a huge integer $N$ with
    $$B_N<(1+\epsilon/8)^{4^{N-1}} |P_1(1,0)|^{4^{N-1}} \textup{ and }
    \frac{\log8}{4^{N-1}}<\epsilon/8.$$
Set
    $$R:=8(1+\epsilon/8)^{4^{N-1}}|P_1(1,0)|^{4^{N-1}}.$$
Since $|A_N|, |B_N|< R/8$, we can choose a very small $\delta>0$ such that
    $$|p_N(t)|, |q_N(t)|<R/8,$$
for any $t$ with $|t|<\delta$. Recall that $F_{n+1}(1,t)=F_{1,t}(F_n(1,t))/t^2$ for all $n\geq 1$, and
    $$ F_{1,t}(z,w) = ((w^2 - t z^2)^2, 4tzw(z-w)(tz-w)).$$
So
    $$|q_{N+1}(t)|=|4p_N(t)q_N(t)(p_N(t)-tq_N(t))(p_N(t)-q_N(t))|<R^4/8,$$
for all $|t|<\delta$, and similarly, $|p_{N+1}(t)|<R^4/8$. Inductively, for all $i\geq 0$ and all $t$ with $|t|<\delta$,  we find
    $$|p_{N+i}(t)|, |q_{N+i}(t)|<R^{4^i}/8.$$
Consequently, as the integer $N$ satisfies $\frac{\log8}{4^{N-1}}<\epsilon/8$, for all $t$ with $|t|<\delta$ and all $n\geq N$,
\begin{eqnarray*}    \frac{1}{{\rm deg}(F_n)}\log\|F_n(1,t)\|&\leq& \frac{\log \max (|tq_n(t)|, |p_n(t)|)}{d\cdot 4^{n-1}}\\
              &\leq& \frac{\log R^{4^{n-N}}-\log 8}{d\cdot4^{n-1}}\\
              &\leq& \frac{\log 8}{d\cdot 4^{N-1}}+\frac{\log |P_1(1,0)|+\log(1+\epsilon/8)}{d}\\
              &<& \frac{\log |P_1(1,0)|}{d}+\epsilon.
\end{eqnarray*}
\qed

\medskip
For the lower bound, we need a preliminary lemma:

\begin{lemma}\label{Lattes t=0 c estimate}
There is a constant $0<c<1$ such that
    $$\frac{\|F_{1,t}(z, w)\|}{\|(z,w)\|^4}\geq c|t|^2$$
for all $t$ with $|t|<1/16$ and all $(z,w )\neq (0,0)$. Consequently,
    $$\frac{1}{{\rm deg}(F_{n+i})}\log\|F_{n+i}(1,t)\|-\frac{1}{{\rm deg}(F_n)}\log\|F_n(1,t)\|\geq \frac{\log c}{4^{n-1}d},$$
for any $i\geq 0$.
\end{lemma}
\proof The homogenous function $F_{1,t} $ is given by
    $$ F_{1,t}(z,w) = ((t z^2 - w^2)^2, 4tzw(z-w)(tz-w)).$$
First assume that $w=1$ and $|z|\leq 1$, so $\|(z,w)\|=1$. Then $\|F_{1,t}(z,1)\|\geq |tz^2-1|^2\geq (1-|t|)^2\geq 1/4>|t|^2$ as $|t|<1/16$. Second,  assume that $z=1$ and $|w|\leq 1$. If $|t-w^2|^2< |t|^2/16$, then one of $|w\pm \sqrt{t}|$ must be $<\sqrt{|t|}/2$. Consequently,
\begin{eqnarray*} \|F_{1,t}(1,w)\|&\geq&|4tw(1-w)(t-w)|\\
                      &\geq&4|t|\left(\sqrt{|t|}-\sqrt{|t|}/2\right)\left(1-\sqrt{|t|}/2-\sqrt{|t|}\right)\left(\sqrt{|t|}-\sqrt{|t|}/2-|t|\right)\\
                      &\geq& |t|^2/16, \textup{ as $|t|<1/16$.}
\end{eqnarray*}
The homogeneity of $F_{1,t}$ allows us to conclude that $\|F_{1,t}(z,w)\|/\|(z,w)\|^4\geq c|t|^2$ for all $(z,w) \not= (0,0)$ and $|t| < 1/16$, with $c = 1/16$.

Noting that $F_{n+1}(1,t)=F_{1,t}(F_n(1,t))/t^2$ for all $n$, we see that
\begin{eqnarray*} \frac{\|F_{n+i}(1,t)\|}{\|F_n(1,t)\|^{4^i}}&=&\left(\frac{\|F_{1,t}(F_{n+i-1})\|}{|t|^2\|F_{n+i-1}\|^{4}}\right)
\left(\frac{\|F_{1,t}(F_{n+i-2})\|}{|t|^2\|F_{n+i-2}\|^{4}}\right)^4\cdots \left(\frac{\|F_{1,t}(F_{n})\|}{|t|^2\|F_{n}\|^{4}}\right)^{4^{i-1}}\\
                                         &\geq& c^{4^i}.
\end{eqnarray*}
Taking logarithms gives the second inequality of this lemma. \qed

\medskip
We use Lemma \ref{Lattes t=0 c estimate} together with the estimates of Lemma \ref{t=0 coefficients growth} to show the corresponding lower bound on the size of $\frac{1}{{\rm deg}(F_n)}\log\|F_n(1,t)\|$ for small $t$.

\begin{lemma}  \label{lower bound first case}
For any given very small $\epsilon>0$, there exists a $\delta>0$ and an integer $N>0$ so that
   $$\frac{1}{{\rm deg}(F_n)}\log\|F_n(1, t)\|-\frac{\log |P_1(1,0)|}{d}>-\epsilon$$
for all $|t|<\delta$ and all $n\geq N$.
\end{lemma}

\proof In contrast with the proof of Lemma \ref{Lattes t=0 upper bound}, we define polynomials $p_n(t)$ and $q_n(t)$ by
    $$F_n(1,t)=(A_n+tp_n(t), tq_n(t)).$$
First, we fix a huge $N>0$ such that
    $$\frac{\log c}{4^{N-1}d}>-\epsilon/8,$$
where $c$ is defined in Lemma \ref{Lattes t=0 c estimate}.
By Lemma \ref{t=0 coefficients growth}, $A_N=P_1(1,0)^{4^{N-1}}$. We pick some very small $\delta>0$ so that for all $t$ with $|t|<\delta$, we have
\begin{eqnarray*} \|F_N(1,t)\|&\geq& |A_N+tp_N(t)|\\
                  &\geq& \left(|P_1(1,0)|(1-\epsilon/8)\right)^{4^{N-1}}.
\end{eqnarray*}
That is,
    $$\frac{1}{{\rm deg}(F_N)}\log\|F_N(1,t)\|\geq \frac{\log|P_1(1,0)|}{d}+\frac{\log(1-\epsilon/8)}{d}.$$
Consequently, by Lemma \ref{Lattes t=0 c estimate}, for any $n\geq N$,
\begin{eqnarray*} \frac{\log\|F_n(1,t)\|}{{\rm deg}(F_n)}&\geq & \left(\frac{\log\|F_n(1,t)\|}{{\rm deg}(F_n)}-\frac{\log\|F_N(1,t)\|}{{\rm deg}(F_N)}\right)+ \frac{\log\|F_N(1,t)\|}{{\rm deg}(F_N)}\\
               &\geq& \frac{\log c}{4^{N-1}d}+\frac{\log\|F_N(1,t)\|}{{\rm deg}(F_N)}\\
               &\geq& -\epsilon/8 +\frac{\log|P_1(1,0)|}{d}+\frac{\log(1-\epsilon/8)}{d}\\
               &>&\frac{\log|P_1(1,0)|}{d}-\epsilon, \textup{ as $\epsilon$ is small.}
\end{eqnarray*}
\qed

\subsection {Convergence near $t_1=0$} \label{t1}
To show the sequence $\frac{1}{{\rm deg}(F_n)}\log\|F_n(t,1)\|$ converges locally uniformly near $t_1=0$,  it suffices to show that for any sufficiently small $\epsilon>0$, there are small $\delta_{\epsilon}>0$ and  integer $N_{\epsilon}$ such that for any $n\geq N_{\epsilon}$ and $|t|<\delta_{\epsilon}$ we have
\begin{equation}\label{lattes Fn 1, 0}-\epsilon<\frac{1}{{\rm deg}(F_n)}\log\|F_n(t,1)\|-\frac{\log |P_1(0,1)|}{d}<\epsilon.
\end{equation}
The proof is similar to the proof of convergence near $t_2=0$ (in \S\ref{near infinity}), but in this case, the lower bound requires more care.

Now we study the coefficients of $F_n(t,1)$. Write
   $$F_n(t,1)=(E_n+D_nt+O(t^2), C_n+O(t)).$$
For the coefficients $C_n$ and $D_n$ and $E_n$, we have the following lemma:
\begin{lemma}\label{s=0 coefficients growth}
For all $n\geq 1$, it has $E_n=P_1(0,1)^{4^{n-1}}$. And for the $C_n$ and $D_n$,
   $$\lim_{n\to \infty} \sup |C_n|^{1/4^{n-1}}, \lim_{n\to \infty} \sup |D_n|^{1/4^{n-1}}\leq |P_1(0,1)|.$$
\end{lemma}

\proof
From the inductive formula $F_{n+1}(t,1)=F_{t,1}(F_n(t,1))$ with
    $$ F_{t,1}(z,w) = ((z^2 - tw^2)^2,4zw(z-w)(z-tw)),$$
we obtain $E_{n+1}=E_n^4$ and $E_1=P_1(0,1)$. So for any $n\geq 1$, $E_n=P_1(0,1)^{4^{n-1}}$.

Let $b_1=0$ and $b_n=2b_{n-1}+1$.  Then inductively, we can write
  $$\begin{cases}C_n=:P_1(0,1)^{4^{n-1}-b_n-1}C^*_n,\\
               D_n=:P_1(0,1)^{4^{n-1}-b_n-1}D^*_n,
   \end{cases}$$
with
 $$\begin{cases}C^*_{n+1}=4C_{n}^*(P_1(0,1)^{b_n+1}- C^*_{n}),\\
               D^*_{n+1}=4P_1(0,1)^{b_n+1}D^*_{n}-2{C^*_{n}}^2,
   \end{cases}$$
for all $n\geq 1$. As both $C_n^*$ and $D_n^*$ grow at most quadratically when $n$ increases by one, we obtain that $\lim_{n\to \infty} \sup |C_n^*|^{1/4^{n-1}}, \lim_{n\to \infty} \sup |D_n^*|^{1/4^{n-1}}\leq 1$. Consequently, we have $\lim_{n\to \infty}\sup |C_n|^{1/4^{n-1}}, \lim_{n\to \infty}\sup |D_n|^{1/4^{n-1}}\leq |P_1(0,1)|$, and the lemma is proved.\qed

\medskip
The growth of the coefficients of $C_n, D_n$ and $E_n$ in Lemma \ref{s=0 coefficients growth} provides a uniform upper bound on the size of $\frac{1}{{\rm deg}(F_n)}\log\|F_n(t,1)\|$ for small $t$.  The proof is exactly the same as the proof of Lemma \ref{Lattes t=0 upper bound}, and we omit the details.

\begin{lemma}\label{Lattes s=0 upper bound}
For any given $\epsilon>0$, there exists a $\delta>0$ and an integer $N>0$ so that
   $$\frac{1}{{\rm deg}(F_n)}\log\|F_n(t,1)\|-\frac{\log |P_1(0,1)|}{d}<\epsilon$$
for all $|t|<\delta$ and all $n\geq N$.
\end{lemma}

\medskip
This next lemma is similar to Lemma \ref{Lattes t=0 c estimate}, and we omit the proof.

\begin{lemma}\label{Lattes s=0 c estimate}
There is a constant $0<c<1$ such that
    $$\frac{\|F_{t,1}(z, w)\|}{\|(z,w)\|^4}\geq c|t|^2$$
for all $t$ with $|t|<1/16$ and all $(z,w )\neq (0,0)$. Consequently,
    $$\frac{1}{{\rm deg}(F_{n+i})}\log\|F_{n+i}(t,1)\|-\frac{1}{{\rm deg}(F_n)}\log\|F_n(t,1)\|\geq \frac{\log (c|t|^2)}{4^{n-1}d},$$
for any $i\geq 0$.
\end{lemma}

\medskip
Obtaining the lower bound on $\frac{1}{{\rm deg}(F_n)}\log\|F_n(t,1)\|$ for small $t$ is  more delicate than in Lemma \ref{lower bound first case}; we use Lemma \ref{Lattes s=0 c estimate} together with the estimates of Lemma \ref{s=0 coefficients growth}.

\begin{lemma}
For any given very small $\epsilon>0$, there exists a $\delta>0$ and an integer $N>0$ so that
   $$\frac{1}{{\rm deg}(F_n)}\log\|F_n(t, 1)\|-\frac{\log |P_1(0,1)|}{d}>-\epsilon$$
for all $|t|<\delta$ and all $n\geq N$.
\end{lemma}

\proof  We define polynomials $p_n(t)$ and $q_n(t)$ by
   $$F_n(t,1)=:(E_n+tp_n(t), q_n(t)).$$
By Lemma \ref{s=0 coefficients growth}, there is a huge integer $N$, such that
   $$|C_N|, |D_N|<(1+\epsilon/16)^{4^{N-1}}|P_1(0,1)|^{4^{N-1}}. $$
 By increasing $N$ if necessary, we may also assume that
   $$\frac{\log(c\epsilon^2/8^2)}{4^{N-1}d}>-\epsilon/10 \textup{ and } \frac{2\log 8}{4^{N-2}d}<\epsilon/8$$
where the constant $c$ is defined in Lemma \ref{Lattes s=0 c estimate}.

Set
   $$R:=8(1+\epsilon/16)^{4^{N-1}}|P_1(0,1)|^{4^{N-1}}.$$
Since $|C_N|, |D_N|<R/8$, we can choose a very small $\delta>0$ such that
   $$|p_N(t)|, |q_N(t)|<R/8$$
for all $t$ with $|t|<\delta.$ Recall that $F_{n+1}(t,1)=F_{t,1}(F_n(t,1))$ for all $n$ and
   $$F_{t,1}(z,w) = ((z^2 - tw^2)^2,4zw(z-w)(z-tw)).$$
So one has
\begin{eqnarray*} |q_{N+1}(t)|&=&|4q_N(t)(E_N+tp_N(t))(E_N+tp_N(t)-q_N(t))\cdot\\
              &&(E_N+tp_N(t)-tq_N(t))|\\
            &<& R^4/8,
\end{eqnarray*}
for all $t$ with $|t|<\delta$. Similarly, we have $|p_{N+1}(t)|<R^4/8$. Inductively, we obtain
\begin{equation}\label{|p_{N+i}(t)|, |q_{N+i}(t)|<R^{4^i}/8}
 |p_{N+i}(t)|, |q_{N+i}(t)|<R^{4^i}/8,
\end{equation}
for all $i\geq 0$ and all $t$ with $|t|<\delta$.

Choose an integer $N'>N$ so that
    $$\delta':=\frac{\epsilon}{8(1+\epsilon/16)^{4^{N'-1}}8^{4^{N'-N}}}<\delta.$$
For any $n\geq N'$ and $t$ with
    $$|t|\leq \frac{\epsilon}{8(1+\epsilon/16)^{4^{n-1}}8^{4^{n-N}}}\leq\delta',$$
as $F_n(t,1)=(E_n+tp_n(t), q_n(t))$, we have
\begin{eqnarray}\label{n's equality} \frac{\log\|F_n(t,1)\|}{{\rm deg}(F_n)}&\geq &\frac{\log|E_n+tp_n(t)|}{d\cdot 4^{n-1}}\\
           &\geq& \frac{\log|P_1(0,1)|}{d}+\frac{\log(1-\epsilon/16)}{d\cdot 4^{n-1}}\textup{ by (\ref{|p_{N+i}(t)|, |q_{N+i}(t)|<R^{4^i}/8})}\\
           &>&\frac{\log|P_1(0,1)|}{d}-\epsilon, \textup{ as $\epsilon$ is small.}  \end{eqnarray}
For any $n\geq N'$ and $t$ with
    $$ \frac{\epsilon}{8(1+\epsilon/16)^{4^{n-1}}8^{4^{n-N}}}<|t|\leq \delta',$$
there is a $j$ with $N'\leq j<n$, such that
$$ \frac{\epsilon}{8(1+\epsilon/16)^{4^{j}}8^{4^{j+1-N}}}<|t|\leq \frac{\epsilon}{8(1+\epsilon/16)^{4^{j-1}}8^{4^{j-N}}}.$$
By the previous argument, as $\frac{\log(c\epsilon^2/8^2)}{4^{N-1}d}>-\epsilon/10,$
    \begin{eqnarray*} \frac{\log\|F_n(t,1)\|}{{\rm deg}(F_n)}&\geq & \left(\frac{\log\|F_n(t,1)\|}{{\rm deg}(F_n)}-\frac{\log\|F_{j}(t,1)\|}{{\rm deg}(F_{j})}\right)+ \frac{\log\|F_{j}(t,1)\|}{{\rm deg}(F_{j})}\\
               &\geq& \frac{\log (c|t|^2)}{4^{j-1}d}+\frac{\log\|F_j(t,1)\|}{{\rm deg}(F_j)}\textup{, by Lemma \ref{Lattes s=0 c estimate}}\\
               &\geq& \left(\frac{\log( c\epsilon^2/8^2)}{4^{j-1}d}-\frac{8}{d}\log(1+\epsilon/16)-\frac{2\log 8}{4^{N-2}d}\right)+\\
               &&\left(\frac{\log|P_1(0,1)|}{d}+\frac{\log(1-\epsilon/16)}{ 4^{j-1}d}\right) \textup{ by (\ref{n's equality}})\\
               &>&\frac{\log|P_1(0,1)|}{d}-\epsilon, \textup{ as $\epsilon$ is small and $N$ is huge.}
\end{eqnarray*}

\qed

\medskip
\subsection{Convergence near $t_1=t_2$} After we set $t_2=1-t$ and $t_1=1$, the proof of the locally uniform convergence statement near $t_1=t_2$ becomes the same as the proof for $t_1=0$ case in \S\ref{t1}.  The role of $P_1(0,1)$ is played by the difference $P_1(1,1)-Q_1(1,1)$, exactly as in the computations for Proposition \ref{res-lat}.  We leave the details for the reader.

This completes the proofs of Theorem \ref{convergence} and Proposition \ref{values of G_C}.

\bigskip
\section{The bifurcation measure}

In this section, we prove Theorem \ref{integral formula} and Proposition \ref{distinct measures}.

\subsection{The measure of maximal entropy}  \label{mu_t}
It is worth observing first that, for any $t\in \Lambda=\mathbb{P}^1-\{0,1,\infty\}$, the measure of maximal entropy of $f_t$ is given by
	$$\mu_t=\frac{2|t(t-1)|\rho_\Lambda(t)}{|z(z-1)(z-t)|}|dz|^2,$$
where $\rho_\Lambda(t)$ is the density of the hyperbolic metric in $\Lambda$, $|dz|^2$ is the Euclidean area element.  Indeed, the fact that the measure of maximal entropy of $f_t$ takes the form $\frac{\phi(t)}{|z(z-1)(z-t)|}|dz|^2$ is well-known; see for example \cite[Appendix A]{Milnor:mating}.  The expression for $\phi(t)$ is due to McMullen, who proved the following identity
$$\frac{1}{\rho_\Lambda(t)}=2|t(t-1)|\int_{\mathbb{P}^1}\frac{|dz|^2}{|z(z-1)(z-t)|}$$
by means of Teichm\"uller theory \cite[Theorem 4.13]{McMullen:course}.

\subsection{Proof of Theorem \ref{integral formula}}
We first consider the constant marked point case. Let $C=(z,w)$ be the lift of the constant marked point $c$, with $z/w\in \mathbb{P}^1-\{0,1,\infty\}$.
In this case, ${\rm gcd}(F_{t_1,t_2}(C))=1$ and  $F_1=F_{t_1,t_2}(C)$ has degree two.  Set $F_{n+1}=F_{t_1,t_2}(F_n)/t_2^2$ for $n\geq1$.  One may verify that $F_n$ is homogenous in $(z,w)$ of degree $4^n$ and homogenous in $(t_1,t_2)$ of degree ${\rm deg}(F_n)=2\cdot4^{n-1}$.  By Theorem \ref{convergence}, the sequence
	$$\frac{1}{{\rm deg}(F_n)}\log\|F_n\|$$
converges locally uniformly in $\mathbb{C}^2-\{(0,0)\}$ to a continuous function, denoted by $G_{(z,w)}(t_1,t_2)$.

For each fixed $t\in  \mathbb{P}^1-\{0,1,\infty\}$, the function $G_{(z,w)}(t,1)$ satisfies
	$$G_{(z,w)}(t,1)=2\log\|(z,w)\|+O(1)$$
when $\|(z,w)\|$ is large.  As a function of $(z,w) \in \C^2\setminus\{0,0\}$, it is  continuous and plurisubharmonic, and it satisfies
	 $$\frac{1}{2}d_{(z,w)}d^c_{(z,w)}G_{(z,w)}(t,1)=\pi^* \mu_t$$
where $\pi: \C^2\setminus\{0,0\} \to \P^1$ is the projection $\pi(z,w) = z/w$ and $\mu_t$ is the measure of maximal entropy for $f_t$; see e.g. \cite[\S4]{Hubbard:Papadopol}.   Writing $G=G_{(z,1)}(t,1)$, the above relation becomes
	 $$\frac{i}{\pi}\partial_z \partial_{\bar{z}} G=2 \mu_t.$$
Therefore, by the discussion in \S\ref{mu_t}, the function $G$ satisfies the Poisson equation:
	$$\frac{1}{\pi}\frac{\partial^2G}{\partial z\partial \bar{z}}=\frac{2|t(t-1)|\rho_\Lambda(t)}{|z(z-1)(z-t)|}.$$
One solution of this equation is
 	 $$I(z,t)=D(t)\int_{\mathbb{P}^1}\frac{\log|z-\zeta|}{|\zeta(\zeta-1)(\zeta-t)|}|d\zeta|^2.$$
where $D(t)=4|t(t-1)|\rho_\Lambda(t)$.

For fixed $t\in \mathbb{P}^1-\{0,1,\infty\}$, note that $z\mapsto I(z,t)$ is continuous on all of $\C$; to handle the singularities of the integrand, we need only observe that $(\log|\zeta|)/|\zeta|$ is integrable near the origin in $\C$.  We have already remarked that $z\mapsto G_{(z,1)}(t,1)$ is continuous at all $z\in \C$.  Now consider the difference $G_{(z,1)}(t,1)-I(z,t)$.
The function $z\mapsto G_{(z,1)}(t,1)-I(z,t)$ is harmonic in $\mathbb{C}-\{0,1,t\}$.
By the removable singularity theorem for harmonic functions, $z\mapsto G_{(z,1)}(t,1)-I(z,t)$ is harmonic in $\mathbb{C}$.
Moreover, one may verify that $G_{(z,1)}(t,1)=O(\log|z|)$ as $|z|\to\infty$. This implies that $G_{(z,1)}(t,1)-I(z,t)$ is a function depending only on $t$, so we may write
   $$G_{(z,1)}(t,1)-I(z,t)=\kappa(t).$$

We can determine $\kappa(t)$ in the following way. Letting $z=0$,  we see that
   $$G_{(0,1)}(t,1)=I(0,t)+\kappa(t)=D(t)\int_{\mathbb{P}^1}\frac{\log|\zeta|}{|\zeta(\zeta-1)(\zeta-t)|}|d\zeta|^2+\kappa(t).$$
It's interesting to note that   $G_{(0,1)}(t,1)$ can be evaluated easily from its definition, as
  $$G_{(0,1)}(t,1)=\log|t|.$$
By changing variables (let $\zeta=t/\xi$), we have that
  $$\int_{\mathbb{P}^1}\frac{\log|\zeta|}{|\zeta(\zeta-1)(\zeta-t)|}|d\zeta|^2=\int_{\mathbb{P}^1}\frac{\log|t|-\log|\xi|}{|\xi(\xi-1)(\xi-t)|}|d\xi|^2.$$
As the total mass of the maximal entropy measure $\mu_t$ for $f_t$ is one, we see that
  	$$I(0,t) =\frac{D(t)}{2}\int_{\mathbb{P}^1}\frac{\log|t|}{|\xi(\xi-1)(\xi-t)|}|d\xi|^2=\log |t|\int_{\mathbb{P}^1}1 d\mu_t=\log|t|.$$
In this way, we find
  	$$\kappa(t)=\log|t|-I(0,t)=0.$$
So we have the following explicit formula
\begin{equation} \label{constant point formula}
  	 G_{(z,1)}(t,1)=D(t)\int_{\mathbb{P}^1}\frac{\log|z-\zeta|}{|\zeta(\zeta-1)(\zeta-t)|}|d\zeta|^2.
\end{equation}
Note that for $\alpha,\beta \in \mathbb{C}^*$,
 	$$G_{(\alpha z,\alpha w)}(\beta t_1, \beta t_2)=G_{(z,w)}(t_1, t_2)+2\log|\alpha|+\log|\beta|.$$
By this scaling property, we have
\begin{eqnarray*} G_{(z,w)}(t_1,t_2)&=&G_{(z/w,1)}(t_1/t_2,1)+2\log|w|+\log|t_2|\\
 	 &=&|t_2| \, D(t_1/t_2)\int_{\mathbb{P}^1}\frac{\log|z-w\zeta|}{|\zeta(\zeta-1)(t_2\zeta-t_1)|}|d\zeta|^2
 			+\log|t_2|.\end{eqnarray*}
providing the formula in homogeneous coordinates for the case of constant starting point $C = (z,w)$.

Now consider the general marked point $C=(c_1(t_1, t_2),c_2(t_1, t_2))$, where $c_1, c_2$ are homogenous polynomials in $(t_1,t_2)$ of the same degree. Comparing the definitions of $G_C(t_1,t_2)$ and $G_{(z,w)}(t_1, t_2)$, we have
 $$G_C(t_1,t_2)=\frac{2}{d}G_{(z,w)}(t_1, t_2)\Big|_{(z,w)=(c_1,c_2)}-\frac{1}{d}\log|{\rm gcd}(F_{t_1,t_2}(C))|.$$
So $G_C$ has the following explicit formula:
\begin{eqnarray*}
 G_C(t_1,t_2)&=&\frac{2}{d}\left(|t_2|D(t_1/t_2) \cdot \int_{\mathbb{P}^1}\frac{\log|c_1-c_2\zeta|}{|\zeta(\zeta-1)(t_2\zeta-t_1)|}|d\zeta|^2+\log|t_2|\right)\\
&& \qquad -\frac{1}{d}\log|{\rm gcd}(F_{t_1,t_2}(C))|.
 \end{eqnarray*}

\subsection{Proof of Proposition \ref{distinct measures}}
Fix point $a \in \C\setminus\{0,1\}$, and set $C = (a,1)$.  Then, as in Theorem \ref{convergence}, we have
	$$F_1 = (P_1(t_1, t_2), Q_1(t_1, t_2)) =  ((t_1 - t_2 a^2)^2, 4t_2a(1-a)(t_1-t_2a))$$
with $d = \deg(F_1) = 2$.  From Theorem \ref{integral formula}, or from equation (\ref{constant point formula}) in its proof, the potential function for the bifurcation measure $\mu_a$ is given explicitly by
	 $$G_a(t,1)=D(t)\int_{\mathbb{P}^1}\frac{\log|a-\zeta|}{|\zeta(\zeta-1)(\zeta-t)|}|d\zeta|^2.$$
From Proposition \ref{values of G_C}, we know that
	$$G_a(1,0)=0,  \quad G_a(0,1)= 2 \log |a|, \quad G_a(1,1)=2\log|1-a|.$$

Now assume we have two points $a, b\in \C\setminus\{0,1\}$ for which $\mu_a = \mu_b$.  Then the harmonic function
	$$G_a(t,1)-G_b(t,1)= G_a(1,0)-G_b(1,0)=0$$
for all $t\in\C$, because $G_a$ and $G_b$ grow at most $O(\log|t|)$ when $t$ is near $\infty$.  This implies
	$$|a|=|b| \quad\mbox{ and } \quad |1-a|=|1-b|.$$
This happens if and only if $a=b$ or $a=\bar{b}$.  However, if we assume $\Im a \not=0$, then
\begin{eqnarray*} &&G_a(t,1)-G_{\bar{a}}(t,1)\\
&=&D(t)\int_{\mathbb{P}^1}\log\left|\frac{a-\zeta}{\bar{a}-\zeta}\right|\frac{|d\zeta|^2}{|\zeta(\zeta-1)(\zeta-t)|}  \\
&=&{D(t)}\int_{\Im\zeta>0} \log\left|\frac{a-\zeta}{\bar{a}-\zeta}\right|\left(\frac{1}{|\zeta-t|}-\frac{1}{|\zeta-\bar{t}|}\right)\frac{|d\zeta|^2}{|\zeta(\zeta-1)|}\\
&\neq& 0 \ \ \mbox{ whenever } \Im t \not=0.
\end{eqnarray*}
This shows that $\mu_a = \mu_b$ if and only if $a=b$, completing the proof of Proposition \ref{distinct measures}.

\bigskip
\section{Equidistribution}

In this section we provide the proofs of Theorems \ref{equidistribution} and \ref{synchrony}.  For the statement of the arithmetic equidistribution theorem, we follow the language and presentation of \cite{BD:polyPCF}, which follows the original treatments of \cite{BRbook, FRL:equidistribution}.  All of the work of Sections \ref{growth} and \ref{estimates} goes towards showing the hypotheses of the equidistribution theorem are satisfied.

\subsection{The arithmetic equidistribution theorem}
Let $k$ be a number field and let $\kbar$ denote a fixed algebraic closure of $k$.  Any number field $k$ is equipped with a set $\cM_k$ of pairwise inequivalent nontrivial absolute values, together with a positive integer $N_v$ for each $v \in \cM_k$, such that
\begin{itemize}
\item for each $\alpha \in k^*$, we have $|\alpha|_v = 1$ for all but finitely many $v \in \cM_k$; and
\item  every $\alpha \in k^*$ satisfies the {\em product formula}
\begin{equation} \label{product formula}
\prod_{v \in \cM_k} |\alpha|_v^{N_v} \ = \ 1 \ .
\end{equation}
\end{itemize}
For each $v \in \cM_k$, let $k_v$ be the completion of $k$ at $v$, let $\kvbar$ be an algebraic closure of $k_v$, and let $\CC_v$ denote the completion of $\kvbar$; for each $v\in \cM_k$, we fix an embedding of $\kbar$ into $\C_v$.  We let $\PP^1_{\Berk,v}$ denote the Berkovich projective line over $\CC_v$, which is a canonically defined path-connected compact Hausdorff space containing $\PP^1(\CC_v)$ as a dense subspace.  If $v$ is archimedean, then $\CC_v \cong \CC$ and $\PP^1_{\Berk,v} = \PP^1(\CC)$.

For each $v \in \cM_k$ there is a distribution-valued Laplacian operator $\Delta$ on $\PP^1_{\Berk,v}$.  For example, the function
$\log^+|z|_v$ on $\PP^1(\CC_v)$ extends naturally to a continuous real valued function $\PP^1_{\Berk,v} \backslash \{ \infty \} \to \RR$ and
\[
\Delta \log^+|z|_v = \delta_{\infty} - \lambda_v,
\]
where $\lambda_v$ is the uniform probability measure on the complex unit circle $\{ |z| = 1 \}$ when $v$ is archimedean and
$\lambda_v$ is a point mass at the Gauss point of $\PP^1_{\Berk,v}$ when $v$ is non-archimedean.

A probability measure $\mu_v$ on $\PP^1_{\Berk,v}$ is said to have {\em continuous potentials} if $\mu_v - \lambda_v = \Delta g$
with $g : \PP^1_{\Berk,v} \to \RR$ continuous.  If $\mu$ has continuous potentials then there is a corresponding {\em  Arakelov-Green function} $g_{\mu} : \PP^1_{\Berk,v} \times \PP^1_{\Berk,v} \to \RR \cup \{ +\infty \}$ which is characterized by the differential equation $\Delta_x g_{\mu}(x,y) = \delta_y - \mu$ and the normalization
\begin{equation} \label{normalization}
	\iint g_{\mu}(x,y) d\mu(x) d\mu(y) = 0.
\end{equation}

An {\em adelic measure} on $\PP^1$ (with respect to the field $k$) is a collection $\mu = \{ \mu_v \}_{v \in M_k}$ of probability measures
on $\PP^1_{\Berk,v}$, one for each $v \in \cM_k$, such that
\begin{itemize}
\item  $\mu_v = \lambda_v$ for all but finitely many $v \in \cM_k$; and
\item  $\mu_v$ has continuous potentials for all $v \in \cM_k$.
\end{itemize}
If $\rho,\rho'$ are measures on $\PP^1_{\Berk,v}$, we define the
{\em $\mu_v$-energy} of $\rho$ and $\rho'$ by
$$( \rho, \rho' )_{\mu_v} := \frac{1}{2} \iint_{\PP^1_{\Berk,v} \times \PP^1_{\Berk,v} \backslash {\rm Diag}} g_{\mu_v}(x,y) d\rho(x) d\rho'(y).$$
Fix a finite subset $S$ of $\PP^1(\kbar)$ which is $\Gal(\kbar/k)$-invariant.  For each $v \in \cM_k$, we denote by $[S]_v$ the discrete probability measure on $\PP^1_{\Berk,v}$ supported equally on elements of $S$.  We let $|S|$ denote the cardinality of $S$.   For each such $S$ with $|S|>1$, the {\em canonical height} of $S$ associated to the adelic measure $\mu = \{ \mu_v \}_{v \in M_k}$ is defined by
\begin{equation}  \label{height definition}
\hhat_{{\mu}}(S) :=  \frac{|S|}{|S|-1}\sum_{v \in \cM_k} N_v \cdot ([S]_v,[S]_v)_{\mu_v}.
\end{equation}
The constants $N_v$ are the same as those appearing in the product formula (\ref{product formula}).

\begin{remark}  The definition of $\hhat_{\mu}$ differs slightly from that given in \cite{BD:polyPCF} or \cite{FRL:equidistribution}; the factor of $|S|/(|S|-1)$ is included to match the definition of the N\'eron-Tate height.  See Proposition \ref{same heights}; compare \cite[Lemma 10.27]{BRbook}.  In fact, this definition of $\hhat_\mu$ extends naturally to include points in $k$, and therefore to define a $\Gal(\kbar/k)$-invariant function $\hhat_\mu: \P^1(\kbar)\to \R$, as can be seen from equation (\ref{point definition}) in the proof of Proposition \ref{same heights}.  
\end{remark}

\begin{theorem} \cite{BRbook, FRL:equidistribution}
\label{arithmetic equidistribution}
Let $\hhat_{{\mu}}$ be the canonical height associated to an adelic measure $\mu$.
Let $\{S_n\}_{n\geq 0}$ be a sequence of finite subsets of $\PP^1(\ksep)$ for which $\hhat_{{\mathbb \mu}}(S_n) \to 0$
as $n \to \infty$.  Then $[S_n]_v$ converges weakly to $\mu_v$ on $\PP^1_{\Berk,v}$ as $n \to \infty$ for all $v \in \cM_k$.
\end{theorem}

\subsection{Computing the Arakelov-Green function}
In the application of Theorem \ref{arithmetic equidistribution} to this article, we will construct height functions explicitly in terms of the homogeneous escape-rate function $G_C$ of Theorem \ref{convergence}.

More generally, suppose that
	$$F_n : (\kbar)^2 \to (\kbar)^2$$
is a sequence of homogeneous polynomial maps defined over a number field $k$.  Assume that $\Res(F_n) \not=0$ for all $n$ and that
	$$\lim_{n\to\infty} \frac{1}{\deg(F_n)} \log\|F_n\|_v$$
converges locally uniformly in $(\C_v)^2 \setminus \{(0,0)\}$ at a given place $v$ of $k$, to a continuous, real-valued function $G_{v}: (\C_v)^2 \setminus \{(0,0)\} \to \mathbb{R}$.  Here, $\|(a,b)\|_v = \max\{|a|_v, |b|_v\}$ as in the archimedean case. Then, in fact, $G_v$ determines a continuous potential function for a probability measure $\mu_v$ on $\P^1_{Berk,v}$.  For $x, y\in \P^1(\C_v)$, the Arakelov-Green function for $\mu_v$ is given by
\begin{equation}  \label{explicit g}
g_{\mu_v}(x,y) = -\log|\tilde{x} \wedge \tilde{y}|_v + G_{v}(\tilde{x}) + G_{v}(\tilde{y}) + \log \capacity(K_v),
\end{equation}
for any choice of lifts $\tilde{x}$ of $x$ and $\tilde{y}$ of $y$ to $(\C_v)^2$.  Here $K_v = \{(a,b)\in (\C_v)^2: G_{v}(a,b) \leq 0\}$. The quantity $\log \capacity(K_v)$ is exactly what is needed to normalize $g_{\mu_v}$ according to (\ref{normalization}).

See \cite[\S10.2]{BRbook} for details, where each step is explained in the setting where $F_n$ is the $n$-th iterate of a homogeneous polynomial map $F_1$ with non-zero resultant.   The quantity $R_v$ appearing in \cite[Lemma 10.10]{BRbook} should be replaced by $\capacity(K_v)$.  The capacity $\capacity(K_v)$ is computed in terms of the resultants of $F_n$ exactly as in Theorem \ref{cap=res}; the same proof works also in the non-archimedean case.  Namely,
\begin{equation}  \label{capacity in all places}
 	\capacity(K_v)=\lim_{n\rightarrow\infty} |{\Res}(F_n)|_v^{-1/\deg(F_n)^2}
\end{equation}
for the given place $v$.

\subsection{The N\'eron-Tate height}  \label{NT sum}
Fix $c\in k(t)\setminus\{0,1,t\}$ for a number field $k$.  Recall from the Introduction that  $P=P_c(t)$ denotes a point in $E_t$ with $x$-coordinate $x(P)=c$, for each $t\in \kbar$. We fix a homogenous lift $C$ of $x(P)$, also defined over $k$, and write $C(t_1,t_2)=(X_1(P),X_2(P))$.  The N\'eron-Tate height of $P$ defines a function $\hhat_c: \P^1(\kbar) \to \R_{\geq 0}$; for $t\in k$, it is given by
\begin{eqnarray*}
\hat{h}_c(t) &=& \hat{h}_{E_t}(P)  \; =\;  \frac12 \lim_{n\to\infty} \frac{1}{4^n} h(x(2^nP)) \\
&=& \frac{1}{2 [k:\Q]}\lim_{n\rightarrow\infty}\frac{1}{4^n} \sum_{v\in \mathcal{M}_{k}} N_v \log\left(\max\{|X_1(2^nP)|_v, |X_2(2^nP)|_v\}\right),
\end{eqnarray*}
where $2^nP$ denotes the image of $P$ under the multiplication-by-$2^n$ map on the elliptic curve, and $h$ is the logarithmic Weil height on $\P^1(\Qbar)$; see, e.g., \cite{Silverman:elliptic}.  If $S$ is a finite, $\Gal(\kbar/k)$-invariant subset of $\P^1(\kbar)$, then we set
$$\hat{h}_c(S) := \frac{1}{|S|} \sum_{t\in S} \hat{h}_c(t). $$

\subsection{Proof of Theorem \ref{equidistribution}}
\label{algebraic case}
Let $c$ be as in the theorem, defined over a number field $k$, and let $C(t_1, t_2)$ be a homogeneous lift of $c$, also defined over $k$.

Fix a place $v$ of $k$.  Define the polynomial maps $F_n$ as in Theorem \ref{convergence}, and set
	$$G_{C,v}(t_1, t_2) = \lim_{n\to\infty}\frac{1}{\deg(F_n)} \log\|F_n(t_1, t_2)\|_v$$
for each place $v$ of $k$, with $(t_1, t_2) \in (\C_v)^2\setminus\{(0,0)\}$.  The proof of convergence in the archimedean case shows {\em mutatis mutandis} that the convergence is locally uniform for all places $v$.  A line-by-line analysis of the proof of Theorem \ref{convergence} shows that the proof uses nothing more than the triangle inequality and elementary algebra.  As such, the estimates can only be improved when the usual triangle inequality is replaced by the ultrametric inequality in the case of a non-archimedean absolute value.

\begin{lemma} \label{good-place} For all but finitely many places $v$,  we have that
 $$\|F_n(t_1, t_2)\|_v = \|(t_1, t_2)\|_v^{\deg(F_n)}$$
 for all $n$ and all $(t_1, t_2) \in (\C_v)^2$. In particular,  the escape-rate function $G_{C,v}$ satisfies
	$$G_{C,v}(t_1, t_2) = \log\|(t_1, t_2)\|_v= \log \max\{|t_1|_v, |t_2|_v\}$$
\end{lemma}

\proof
The coefficients of $F_1 = (P_1, Q_1)$ lie in $k$.  Let $v$ be any non-archimedean place for which all coefficients of $P_1$ and $Q_1$ lie in the valuation ring of $k$ (and so have absolute value $\leq 1$) and such that
	$$|2|_v = |\Res(F_1)|_v = |P_1(0,1)|_v = |P_1(1,0)|_v = |P_1(1,1)-Q_1(1,1)|_v = 1.$$
Note that this holds for all but finitely many places $v$.  Then, by Proposition \ref{res-lat}, we have
	$$|\Res(F_n)|_v = 1$$
for all $n$.

Applying \cite[Lemma 10.1]{BRbook}, we may conclude that
	$$\|F_n(t_1, t_2)\|_v = \|(t_1, t_2)\|_v^{\deg(F_n)}.$$
for all $n$ and all $(t_1, t_2) \in (\C_v)^2$.   This completes the proof of the lemma.
\qed

\bigskip
It follows that the functions $\{G_{C,v}\}_{v\in \cM_k}$ extend to the Berkovich spaces $\P^1_{Berk, v}$ to define a collection of continuous (homogeneous) potential functions for an adelic measure $\mu = \{\mu_{C,v}\}$.   We call this measure $\mu$ the {\em adelic bifurcation measure} of $c$.   Using (\ref{capacity in all places}), we obtain a formula for the capacity of $K_v = \{G_{C,v} \leq 0\} \subset (\C_v)^2$, exactly as in Theorem \ref{capacity}, with exactly the same proof:
	$$\capacity(K_v)  = |4|_v^{-\frac{1}{3 d}}\left|\frac{(Q(1,1)-P(1,1))P(0,1)}{P(1,0)^2}\right|_v^{-\frac{1}{3d^2}} |\Res(F_1)|_v^{-\frac{1}{d^2}}$$
where $F_1 = (P, Q)$ and $d =d(c)$ are defined in Theorem \ref{convergence}.  This explicit formula for the capacity and the product formula (\ref{product formula}) imply that
\begin{equation} \label{capacity product}	
	\prod_v  \capacity(K_v)^{N_v} = 1.
\end{equation}

The next lemma shows that the N\'eron-Tate height function $\hhat_c$ may be expressed in terms of the functions $\{G_{C,v}\}$; compare to the definition of the canonical height given in equation (10.40) of \cite{BRbook}. 

\begin{lemma}\label{neron-tate}  The N\'eron-Tate height $\hat{h}_c(t) = \hat{h}_{E_t}(P_c(t))$  satisfies
	 $$\hat{h}_c(t)=\frac{d}{8[k:\Q]}\sum_{v\in \mathcal{M}_k}N_vG_{C,v}(t_1,t_2)$$
for all $t\in k$ and $(t_1, t_2) \in k^2$ with $t_1/t_2 = t$; and
	$$\hat{h}_c(S) = \frac{d}{8[k:\Q]|S|} \sum_{t\in S} \sum_{v\in \mathcal{M}_k}N_vG_{C,v}(t_1,t_2)$$
for all finite, $\Gal(\kbar/k)$-invariant sets $S$,
where $d = d(c)$ is defined in Theorem \ref{convergence}.
\end{lemma}

\proof
By the product formula, the quantity
$$H_k(x)=\prod_{v\in \mathcal{M}_k} \max\{|X_1|_v, |X_2|_v\}^{N_v}$$
does not depend on the lift $\tilde{x}=(X_1,X_2)\in k^2$ of $x\in\P^1(k)$.  As we know that $F_n(t_1,t_2)$ (defined in Theorem \ref{convergence}) is a lift of $x(2^nP)=f_t^{n}(x(P))$, the N\'eron-Tate height $\hat{h}_c(t) = \hat{h}_{E_t}(P)$ is given by
$$\hat{h}_c(t) = \frac{1}{2[k:\Q]}\lim_{n\rightarrow\infty}\frac{1}{4^n}\sum_{v\in \mathcal{M}_k} N_v
\log\|F_n(t_1,t_2)\|_v,$$
for $t\in k$ and $(t_1, t_2) \in k^2$ with $t_1/t_2 = t$, from the expression given in \S\ref{NT sum}.  Now, we decompose $\mathcal{M}_k$  into two sets: $\mathcal{M}_k^1$, consisting of the places  $v$ for which
 $$\|F_n(t_1, t_2)\|_v = \|(t_1, t_2)\|_v^{\deg(F_n)}$$
 for all $n$ and all $(t_1, t_2) \in (\C_v)^2$ (as in Lemma \ref{good-place}),  and $\mathcal{M}_k^2=\mathcal{M}_k\setminus \mathcal{M}_k^1$, which is a finite set. As we have shown that the sequence  $\frac{1}{\deg(F_n)} \log\|F_n\|_v$ with $\deg(F_n)=4^{n-1}d$ converges locally uniformly in $(\C_v)^2 \setminus \{(0,0)\}$ at each place $v$ of $k$, to a continuous, real-valued function $G_{C,v}: (\C_v)^2 \setminus \{(0,0)\} \to \mathbb{R}$, we have
\begin{eqnarray*}
\hat{h}_c(t)
&=&\frac{1}{2[k:\Q]} \lim_{n\to\infty} \frac{1}{4^n}  \left(  \sum_{v\in \mathcal{M}_k^1}N_v \log\|(t_1, t_2)\|_v^{d4^{n-1}} + \sum_{v\in \mathcal{M}_k^2} N_v
\log\|F_n(t_1,t_2)\|_v  \right)  \\
&=& \frac{1}{2[k:\Q]} \lim_{n\to\infty} \left( \frac{d}{4} \sum_{v\in \mathcal{M}_k^1}N_v \log\|(t_1, t_2)\|_v + \frac{1}{4^n} \sum_{v\in \mathcal{M}_k^2} N_v \log\|F_n(t_1,t_2)\|_v  \right)  \\
&=&\frac{d}{8[k:\Q]}\sum_{v\in \mathcal{M}_k}N_vG_{C,v}(t_1,t_2),
\end{eqnarray*}
for $t\in k$ and all $(t_1, t_2) \in k^2$ with $t_1/t_2 = t$.

Now suppose that $S$ is a finite and $\Gal(\kbar/k)$-invariant subset of $\P^1(\kbar)$.  Let $k'$ be a finite extension of $k$ containing $S$.  The proof above shows that
	$$\hhat_c(t) = \frac{d}{8[k':\Q]}\sum_{v\in \mathcal{M}_{k'}}N_vG_{C,v}(t_1,t_2)$$
for all $t\in S$ and $(t_1, t_2)\in (k')^2$ with $t_1/t_2 =t$.  Recall that, for each $v\in \cM_k$, we fixed an embedding of $\kbar$ into $\C_v$, thus making a choice of extension of the absolute value $|\cdot|_v$ to $k'$.  By the Galois invariance of $S$, the sum
	$$\sum_{t\in S} \sum_{v\in \cM_k} N_vG_{C,v}(t_1,t_2)$$
is independent of the choice of embedding for each $v$, and therefore,
\begin{eqnarray*}
 \frac{d}{8|S|} \sum_{t\in S} \frac{1}{[k:\Q]} \sum_{v\in \mathcal{M}_k}N_vG_{C,v}(t_1,t_2)
 &=&  \frac{d}{8|S|} \sum_{t\in S} \frac{1}{[k':\Q]} \sum_{v\in \mathcal{M}_{k'}}N_vG_{C,v}(t_1,t_2)  \\
 &=& \frac{1}{|S|}\sum_{t\in S} \hhat_c(t)
\end{eqnarray*}
\qed

\medskip
We now prove the precise relation between the N\'eron-Tate height $\hhat_c$ and the measure-theoretic height $\hhat_\mu$, defined in (\ref{height definition}), for the adelic bifurcation measure of $c$.  

\begin{prop}\label{same heights}
Let $k$ be a number field; let $S$ be any finite set in $\P^1(\kbar)$ that is $\Gal(\kbar/k)$-invariant with $|S|>1$; and let $\mu = \{\mu_{c,v}\}$ be the adelic bifurcation measure for $c\in k(t) \setminus\{0,1,t\}$.  The canonical height $\hat{h}_\mu$ of (\ref{height definition}) and the  N\'eron-Tate height $\hat{h}_c$ satisfy
	 $$\hhat_{\mu}(S)=\frac{8[k:\Q]}{d} \, \hat{h}_c(S),$$
where $d = d(c)$ is defined in Theorem \ref{convergence}.
\end{prop}

\proof
For each $x\in S$, fix any point $\tilde{x}\in (\kbar)^2$ over $x$.  We may verify directly from the definition of $\hhat_{\mu}$ (exactly as in \cite[Lemma 10.27]{BRbook}):

\begin{eqnarray}
\hhat_{\mu}(S) &=& \frac{|S|}{|S|-1}\sum_{v\in \cM_k} N_v \cdot ([S]_v, [S]_v)_{\mu_v} \nonumber \\
	&=&  \frac{|S|}{|S|-1}  \sum_{v\in \cM_k} \frac{N_v }{2|S|^2}\sum_{x,y\in S, x\not=y}  g_{\mu_v}(x,y) \nonumber \\	
	&=&   \frac{1}{2|S|(|S|-1)} \sum_{ x,y\in S, x\not=y }   \sum_{v\in \cM_k} N_v \cdot ( -\log|\tilde{x} \wedge \tilde{y}|_v + G_{C,v}(\tilde{x}) + G_{C,v}(\tilde{y}) + \log \capacity(K_v) ) \nonumber \\
	&=&  \frac{1}{2|S|(|S|-1)} \sum_{x,y\in S, x\not=y}  \; \sum_{v\in \cM_k} N_v \cdot (G_{C,v}(\tilde{x}) + G_{C,v}(\tilde{y}) ) \label{explain} 
			 \\
	&=& \frac{1}{2|S|(|S|-1)} \sum_{v\in \cM_k} N_v  \cdot 2\, (|S|-1)\sum_{x\in S} G_{C,v}(\tilde{x}) \nonumber \\
	&=&  \frac{1}{|S|} \sum_{v\in \cM_k} N_v \cdot \sum_{x\in S} G_{C,v}(\tilde{x}) \nonumber \\
     &=&  \frac{1}{|S|}  \cdot \sum_{x\in S}\sum_{v\in \cM_k} N_v G_{C,v}(\tilde{x}),   \label{point definition}
\end{eqnarray}
where (\ref{explain}) follows from the product formula and (\ref{capacity product}).  We remark that the formula (\ref{point definition}) shows that $\hhat_\mu$ can be extended naturally to sets $S$ of cardinality 1, and therefore, to define a $\Gal(\kbar/k)$-invariant function on $\P^1(\kbar)$.  By Lemma \ref{neron-tate}, we conclude that
$$\hhat_{\mu}(S)=\frac{8[k:\Q]}{d|S|}\sum_{t\in S}\hat{h}_c(t) = \frac{8[k:\Q]}{d} \; \hat{h}_c(S).$$
\qed

\medskip
With Proposition \ref{same heights}, the conclusion of Theorem \ref{equidistribution} follows immediately from Theorem \ref{arithmetic equidistribution}.

\subsection{Proof of Theorem \ref{synchrony} (and the case $a, b \in \Qbar(t)$ of Theorem \ref{synchrony over C}).} \label{proof of synchrony}
Suppose that $c\in k(t)$, where $k$ is a number field with $c\not= 0, 1, t$.  Recall that
\begin{eqnarray*}
\Tor(c) &=&\{t\in \C\setminus\{0,1\}:  c(t)  \mbox{ is a torsion point on } E_t \}\\
 &=&\{t\in \C\setminus\{0,1\}:  c(t)  \mbox{ is preperidic for } f_t \}.
\end{eqnarray*}
From the convergence result, Theorem \ref{convergence}, we may conclude that $c$ has a nontrivial bifurcation measure $\mu_c$.  It follows that the sequence of functions, $t \mapsto f_t^n(c(t))$, $n\geq 1$, do not form a normal family on $\C\setminus\{0,1\}$; see e.g., \cite[Theorem 9.1]{D:lyap}.  Consequently, with standard arguments in complex dynamics using Montel's theory of normal families, the point $c(t)$ must be preperiodic for $f_t$ for infinitely many $t$; compare \cite[Lemma 2.4]{BD:polyPCF}, \cite[Lemma 2.3]{Dujardin:Favre:critical}.  In other words, $|\Tor(c)| = \infty$.

Now let $a\not=b$ be two such points in $k(t)$.  The argument above explains why condition (2) (that $\Tor(a) = \Tor(b)$) implies condition (1) (that $|\Tor(a)\cap \Tor(b)| = \infty$).

Now we check $(1)\Rightarrow(3)$. Recall  that $P_c(t)$ is torsion on $E_t$ if  and only if $\hat{h}_{E_t}(P_c(t))=0$; see \cite[Theorem 9.3]{Silverman:elliptic}. Thus  for all $t\in \Tor(a)\cap\Tor(b)$, we have $\hat{h}_{a}(t)=0$ and $\hat{h}_{b}(t)=0$. The implication then follows.

To see that $(3)\Rightarrow(2)$, we assume there is a non-repeating infinite sequence $\{t_n\}\subset \kbar$ for which $\hat{h}_{a}(t_n), \hat{h}_{b}(t_n)\rightarrow0$.  Let $S_n=\Gal(\kbar/k)\cdot t_n$.  Since the N\'eron-Tate height is invariant under Galois conjugation, we see that
  $$\hat{h}_a(S_n) \to 0 \qquad \mbox{ and } \qquad \hat{h}_b(S_n) \to 0.$$
From Proposition \ref{same heights}, it follows that
    $$\hat{h}_{\mu_a}(S_n) \to 0 \qquad \mbox{ and } \qquad \hat{h}_{\mu_b}(S_n) \to 0$$
for the canonical heights associated to the adelic bifurcation measures $\{\mu_{a,v}\}$ and $\{\mu_{b,v}\}$.
From the equidistribution of $S_n$ in Theorem \ref{equidistribution}, we deduce that the measures $\mu_{a,v} = \mu_{b,v}$ for all places $v$ of $k$.  It follows that $\hat{h}_{\mu_a} = \hat{h}_{\mu_b}$, and therefore
  $$\hat{h}_{E_t}(P_a(t)) = \hhat_a(t) = \hhat_b(t) = \hhat_{E_t}(P_b(t))$$
for all $t\in \kbar$.  By \cite[Theorem 9.3]{Silverman:elliptic},
 $$\Tor(a)=\{t\in \kbar\setminus\{0,1\}: \hhat_a(t)=0\}=\{t\in \kbar\setminus\{0,1\}: \hhat_b(t)=0\}=\Tor(b).$$

Finally, assume that $a,b\in \Qbar\setminus\{0,1\}$ are constant points.  Assuming conditions (1)--(3), the proof above implies that the archimedean components of the bifurcation measures must coincide.  In particular, $\mu_a = \mu_b$, for the bifurcation measures defined in the Introduction (after Theorem \ref{convergence}).  By Proposition \ref{distinct measures}, we may conclude that $a=b$.  This completes the proof of Theorem \ref{synchrony}.  We have also completed the proof of Theorem \ref{synchrony over C} in the case where $a, b\in \Qbar(t)$.

\bigskip
\section{The proof of Theorem \ref{synchrony over C}}
\label{complex case}

In this section, we complete the proof of Theorem \ref{synchrony over C}.   The result in the case where $a, b\in \Qbar(t)$ was obtained in \S\ref{proof of synchrony}.  The main ideas of this section are inspired by the method of \cite{BD:preperiodic} to pass from results over $\Qbar$ to results over $\C$.

\subsection{The height function}  \label{function field height}
Fix a point $c \in \C(t)$, and assume that $c \not= 0,1, t$.  Assume that $c \not\in \Qbar(t)$ and choose a finitely-generated extension field $k$ over $\Qbar$ so that $c \in k(t)$.  

Viewed as a function field over $\Qbar$, $k$ may be equipped with a product-formula structure; that is, with a collection of absolute values $\cM_k$ satisfying the product formula (\ref{product formula}) for all nonzero elements of $k$.  The ``naive" height function on $\P^1(\kbar)$ is defined by
	$$h(\alpha:\beta) = \frac{1}{[K:k]} \; \sum_{\sigma \in \Gal(K/k)} \; \sum_{v\in \cM_k} N_v \log\|\sigma (\alpha,\beta)\|_v$$
for any finite extension $K$ of $k$ containing $\alpha$ and $\beta$, exactly as in the number field case.  Here, $\|(\alpha, \beta)\|_v = \max\{|\alpha|_v, |\beta|_v\}$.  

Choose any homogeneous lift $C$ of $c$ defined over $k$, as in \S\ref{results}.  From Proposition \ref{persistent}, we know that $c(t)$ is not persistently preperiodic for $f_t$, for $t\in \C\setminus\{0,1\}$.  Theorem \ref{convergence} implies that the potential function $G_C$ is well defined and continuous on $\C^2\setminus\{(0,0)\}$.  

The arguments of \S\ref{algebraic case} carry over to this setting, showing that the escape-rate functions $G_{C,v} : (\C_v)^2 \to \R\cup\{-\infty\}$ are well defined and continuous at all places $v$ of $k$.  In particular, Lemma \ref{neron-tate} applies to show that the height function 
	$$\hat{h}_c(t) :=  \frac{1}{2} \lim_{n\to\infty} \frac{1}{4^n} h(f_t^n(c(t)))$$
on $\P^1(\kbar)\setminus\{0,1,\infty\}$ satisfies
\begin{equation}\label{neron-tate height1}
\hat{h}_c(S) = \frac{d}{8|S|} \sum_{(t_1: t_2) \in S} \sum_{v\in \mathcal{M}_k}N_vG_{C,v}(t_1,t_2)
\end{equation}
for all finite, $\Gal(\kbar/k)$-invariant sets $S$ in $\P^1(\kbar)$, where $d = d(c)$ is the degree defined in Theorem \ref{convergence}.   Consequently, the equidistribution of points of small height, Theorem \ref{arithmetic equidistribution}, holds in this setting as well. 

One key difference between the number field case and this function field setting is that there is no longer an archimedean place of $k$; the bifurcation measure $\mu_c$ on $\P^1(\C)$ is {\em not} a component of the adelic measure $\{\mu_{c,v}\}$.  Another key difference is that 
	$$\Tor(c) := \{t: c(t) \mbox{ is preperiodic for } f_t\} \subset \{t: \hat{h}_c(t)=0\},$$
as can be seen from the definition of $\hat{h}_c$, but we may not have equality of sets.  Indeed, if $t\in \Qbar \subset k$ is such that $c(t)\in \Qbar$, then $\hat{h}_c(t) = 0$ independent of whether $P_c(t)$ is a torsion point on $E_t$.  

\subsection{Proof of Theorem \ref{synchrony over C}}
Fix $a,b \in \C(t)$ as in the statement of the theorem, and let $k$ be a finitely-generated field extension of $\Qbar$ containing the coefficients of $a$ and $b$.  We may assume that at least one of the points, say $a$, is not in $\Qbar(t)$.  Assume that $\Tor(a) = \Tor(b)$.  Exactly as in the beginning of the proof of Theorem \ref{synchrony}, we see that $|\Tor(a)| = \infty$, so that $|\Tor(a)\cap \Tor(b)|=\infty$.  

Now, assuming that $|\Tor(a) \cap \Tor(b)| = \infty$, we may appeal to the main result of \cite{Masser:Zannier:2} to deduce that $P_a$ and $P_b$ are linearly related.  That is, there exist integers $m, n$, not both zero, such that $[m]\cdot P_a + [n]\cdot P_b= 0$ on $E_t$ for all $t\in \C\setminus\{0,1\}$.  In fact, by our assumption that $a, b \not= 0,1,t, \infty$, we see that both $m$ and $n$ must be nonzero (since neither $P_a$ nor $P_b$ can be torsion for all $t$).  

With $P_a$ and $P_b$ linearly related, we immediately deduce that $\Tor(a) = \Tor(b)$, and therefore we have equivalence of conditions (1), (2), and (4).  Assuming condition (1), and since $\Tor(a) \subset \{t: \hhat_a(t)=0\}$ for any $a$, we obtain condition (3) with an infinite sequence of torsion parameters.

Finally, we will show that condition (3) implies condition (1).  We remarked in \S\ref{function field height} that the equidistribution theorem holds in this setting.  Consequently, condition (3) implies that we have equality of measures $\mu_{a,v} = \mu_{b,v}$ at all places $v$ of $k$.  Thus $\hhat_{\mu, a} = \hhat_{\mu,b}$ on $\P^1(\kbar)$.  By (\ref{point definition}) and (\ref{neron-tate height1}), we find that $d(b)\cdot \hat{h}_a = d(a)\cdot \hat{h}_b$, where the degrees $d(a), d(b)$ are defined in Theorem \ref{convergence}.  In particular, $\hhat_a$ and $\hhat_b$ vanish at the same points $t\in \kbar$.

Note that the set of $t \in \Qbar$ such that $a(t) \in \Qbar$ is finite; indeed, if it were infinite, then the algebraic curve parameterized by $(t, a(t))$ in $\C^2$ would be defined over $\Qbar$, contradicting our assumption that $a\not\in \Qbar(t)$.  Note further that for $t\in \Qbar$ with $a(t) \not\in\Qbar$, then $t\not\in \Tor(a)$, since the torsion points of $E_t$ must lie in $\Qbar$.  Therefore, since $|\Tor(a)| = \infty$, there must be infinitely many torsion parameters $t\not\in \Qbar$.

For each $t\not\in \Qbar$, the elliptic curve $E_t/k$ is not isotrivial; that is, it is not isomorphic (over $k$ or even over $\kbar$) to an elliptic curve defined over the constant field $\Qbar$.  (This can be seen explicitly by observing that the $x$-coordinates of the 2-torsion points $0, 1, t, \infty$ on $E_t$ have cross-ratio $t$ that is not in the constant field $\Qbar$.  The cross-ratio, up to reordering of the points, is an isomorphism invariant.)  Therefore, by standard theory of non-isotrivial elliptic curves over function fields \cite[Theorem  III.5.4]{Silverman:Advanced}, or by appealing to the corresponding dynamical statement about $f_t$ \cite[Theorem 1.6]{Baker:functionfields}, we deduce that $\hat{h}_c(t) = 0$ if and only if $t \in \Tor(c)$ whenever $t\not\in \Qbar$, for any choice of $c\in k(t)$.  

Therefore, for each of the (infinitely many) $t\in \Tor(a)$ with $t\not\in \Qbar$, it must be that $t\in \Tor(b)$ since $d(b)\cdot \hat{h}_a =d(a) \cdot \hat{h}_b$.  This concludes the proof that (3) implies (1) and the proof of Theorem \ref{synchrony over C}.

\bigskip \bigskip

\def\cprime{$'$}

\end{document}